\documentclass[11pt]{article}
\usepackage{amsmath, amssymb, epsfig, eepic, bm, amsbsy}
\usepackage{amsfonts}

\newcommand{\bt}   {\begin{theorem}}
\newcommand{\et}   {\end  {theorem}}
\newcommand{\bl}   {\begin{lemma}}
\newcommand{\el}   {\end  {lemma}}
\newcommand{\bp}   {\begin{prop}}
\newcommand{\ep}   {\end  {prop}}
\newcommand{\bc}   {\begin{cor}}
\newcommand{\ec}   {\end  {cor}}
\newcommand{\bd}   {\begin{defn}}
\newcommand{\ed}   {\end  {defn}}

\newcommand{\ba}   {\begin{array}}
\newcommand{\ea}   {\end  {array}}
\newcommand{\be}   {\begin{enumerate}}
\newcommand{\ee}   {\end  {enumerate}}
\newcommand{\bi}   {\begin{itemize}}
\newcommand{\ei}   {\end  {itemize}}

\def\eq#1\en{\begin{equation}#1\end{equation}}
\def\eqsplit#1\ensplit{
    \begin{equation}\begin{split}#1\end{split}\end{equation}
    }
\def\eqalign#1\enalign{
    \begin{align}#1\end{align}
    }
\def\eqmul#1\enmul{
    \begin{multline}#1\end{multline}
    }
\newcommand{\eqarrstar} {\begin{eqnarray*}}
\newcommand{\enarrstar} {\end{eqnarray*}}
\newcommand{\eqarray}   {\begin{eqnarray}}
\newcommand{\enarray}   {\end{eqnarray}}

\newcommand{\lbeq}[1]  {\label{e:#1}}
\newcommand{\refeq}[1] {\eqref{e:#1}}    

\makeatletter
\newcommand{\labelcounter}[2]{{%
    \stepcounter{#1}
    \protected@write\@auxout{}%
    {\string\newlabel{#2}{{\csname the#1\endcsname}{\thepage}}}%
    {\ref{#2}}
    }}
\makeatother
\newcommand{\sss}   { \scriptscriptstyle }

\newcommand{\spose}[1] {{\hbox to 0pt{#1\hss}} }
\newcommand{\ltapprox} {\mathrel{\spose{\lower 3pt\hbox{$\mathchar"218$}}
 \raise 2.0pt\hbox{$\mathchar"13C$}}}
\newcommand{\gtapprox} {\mathrel{\spose{\lower 3pt\hbox{$\mathchar"218$}}
 \raise 2.0pt\hbox{$\mathchar"13E$}}}

\newcommand{\toinf} {\rightarrow \infty}

\newcommand{\bra}[1]    {\left \langle #1 \right |}

\newcommand{\expec}[1]  {\left \langle #1 \right \rangle}

\newcommand{\eps}{\varepsilon}

\newcommand{\nn}{\nonumber}

\numberwithin{equation}{section}

\newcommand{\smallsup}[1] {{\scriptscriptstyle{(\kern-.03em{#1}}\kern-.03em)}}

\newcommand{\indic}[1]{{\bf 1}_{\bra{#1}}} 

\newcommand{\br}[1]{\left( #1 \right)}       
\newcommand{\brh}[1]{\left[ #1 \right]}    
\newcommand{\brc}[3]{\left#1 #3 \right#2}  
\newcommand{\eqa}{\begin{eqnarray}}
\newcommand{\ena}{\end{eqnarray}}
\newcommand{\qed}    {\hfill$\Box$}
\newcommand{\eqan}[1]{\eqalign #1 \enalign}

\newcommand{\vep}{\varepsilon}
\newcommand{\prob}{\mathbb P}

\renewcommand{\expec}{\mathbb E}
\newcommand{\mexpec}[1]{\mathbb E\hspace{-.2em}\brh{#1}}

\newcommand{\eqn}[1]{\begin{equation} #1 \end{equation}}

\newcommand{\cvepone}{C_{\vep}^{\sss(1)}}
\newcommand{\cveptwo}{C_{\vep}^{\sss(2)}}
\newcommand{\cvepthree}{C_{\vep}^{\sss(3)}}
\newcommand{\const}[1]{C_{#1}}

\newtheorem{theorem}{Theorem}[section]
\newtheorem{conj}[theorem]{Conjecture}
\newtheorem{lemma}[theorem]{Lemma}
\newtheorem{prop}[theorem]{Proposition}
\newtheorem{cor}[theorem]{Corollary}

\newtheorem{rem}[theorem]{Remark}

\newcommand{\TW}{{\rm TW}}

\usepackage{color}

\title{A preferential attachment model with random initial degrees}
\author{Maria Deijfen
\thanks{Stockholm University, Department of Mathematics, 106 91 Stockholm,
Sweden. E-mail: {\tt mia@math.su.se}}
\and Henri van den Esker\footnote{Delft University of Technology,
Electrical Engineering, Mathematics and Computer Science, P.O. Box
5031, 2600 GA Delft, The Netherlands. E-mail: {\tt
H.vandenEsker@ewi.tudelft.nl}, {\tt G.Hooghiemstra@ewi.tudelft.nl}
} \and Remco van der Hofstad\footnote{Department of Mathematics
and Computer Science, Eindhoven University of Technology, P.O.\
Box 513, 5600 MB Eindhoven, The Netherlands. E-mail: {\tt
rhofstad@win.tue.nl}} \and Gerard Hooghiemstra$^\dagger$}

\date{June 2020}

\begin{document}

\maketitle

\begin{abstract}
\noindent In this paper, a random graph process $\{G(t)\}_{t\geq
1}$ is studied and its degree sequence is analyzed. Let
$\{W_t\}_{t\geq 1}$ be an i.i.d.\ sequence. The graph process is
defined so that, at each integer time $t$, a new vertex, with
$W_t$ edges attached to it, is added to the graph. The new edges
added at time $t$ are then preferentially connected to older
vertices, i.e., conditionally on $G(t-1)$, the probability that a
given edge is connected to vertex $i$ is proportional to
$d_i(t-1)+\delta$, where $d_i(t-1)$ is the degree of vertex $i$ at
time $t-1$, independently of the other edges. The main result is
that the asymptotical degree sequence for this process is a power
law with exponent $\tau=\min\{\tau_{\sss {\rm W}}, \tau_{\sss {\rm
P}}\}$, where $\tau_{\sss {\rm W}}$ is the power-law exponent of
the initial degrees $\{W_t\}_{t\geq 1}$ and $\tau_{\sss {\rm P}}$
the exponent predicted by pure preferential attachment. This
result extends previous work by Cooper and Frieze, which is
surveyed.
\end{abstract}

\section{Introduction}\label{sec:rules}

Empirical studies on real life networks, such as the Internet, the
World-Wide Web, social networks, and various types of
technological and biological networks, show fascinating similarities.
Many of the networks are {\it small worlds}, meaning that typical
distances in the network are small, and many of them have {\it
power-law degree sequences}, meaning that the number of vertices
with degree $k$ falls off as $k^{-\tau}$ for some exponent
$\tau>1$. See \cite{FFF99} for an example of these phenomena in
the Internet, and \cite{KRRSTU00, KRRT99} for an example on the
World-Wide Web. Also, Table 3.1 in \cite{NewMan} gives an overview
of a large number of networks and their properties.

Incited by these empirical findings, random graphs have been
proposed to model and/or explain these phenomena -- see
\cite{Durr06} for an introduction to random graph models for
complex networks. Two particular classes of models that have been
studied from a mathematical viewpoint are (i) graphs where the
edge probabilities depend on certain weights associated with the
vertices, see e.g.\ \cite{BolJanRio06,BDM-L05,CL1,CL2,S}, and (ii)
so-called preferential attachment models, see e.g.\
\cite{barabasi-1999-286,BBCR03,BR03, BRST01,cf}. The first class
can be viewed as generalizations of the classical
Erd\H{o}s-R\'enyi graph allowing for power-law degrees. In
\cite{BDM-L05}, for instance, a model is considered in which each
vertex $i$ is assigned a random weight $W_i$ and an edge is drawn
between two vertices $i$ and $j$ with a probability depending on
$W_iW_j$. This model, which is referred to as the {\it generalized
random graph}, leads to a graph where vertex $i$ has an asymptotic
degree distribution equal to a Poisson random variable with
(random) parameter $W_i$ as the number of vertices tends to
infinity, that is, the asymptotic degree of a vertex is determined
by its weight. We refer to \cite{Boll01, JanLucRuc00} for
introductions to classical random graphs.

Preferential attachment models are different in spirit in that
they are dynamic, more precisely, a new vertex is added to the
graph at each integer time. Each new vertex comes with a number of
edges attached to it and these edges are connected to the old
vertices in such a way that vertices with high degree are more
likely to be attached to. It can be shown that this leads to graphs
with power-law degree sequences. Note that, in preferential
attachment models, the degree of a vertex increases over time,
implying that the oldest vertices tend to have the largest
degrees.

Consider the degree of a vertex as an indication of its
{\it success}, so that vertices with large degree
correspond to successful vertices. In preferential attachment
models, vertices with large degrees are the most likely vertices
to obtain even larger degrees, that is, successful vertices are
likely to become even more successful. In the literature this is
sometimes called the {\it rich-get-richer} effect. In the
generalized random graph, on the other hand, a vertex is {\it
born} with a certain weight and this weight determines the degree
of the vertex, as described above. This will be referred to as the
{\it rich-by-birth} effect in what follows.

Naturally, in reality, both the rich-get-richer and the
rich-by-birth effect may play a role. To see this, consider for
instance a social network, where we identify vertices with
individuals and edges with social links, that is, an edge is added
between two individuals if they have some kind of social relation
with each other. Then, indeed, we would expect to see both
effects: Firstly, the rich-get-richer effect should be apparent,
since individuals with a high number of social links will in time
acquire more new social links than individuals with few social
links. Evidently, more social contacts imply that the individual
is socially more active, so that he/she meets more people, and, in
turn, each meeting offers a possibility to create a lasting social
link. Thus, these individuals are more likely to get to acquainted
with even more individuals. Secondly, the rich-by-birth effect
comes in due to the fact that some individuals are better in
turning a meeting into a lasting social link than others. The
social activity and skill varies from individual to individual and
could be measured for instance by weights associated with the
individuals. Hence, in reality, both the previous success of a
vertex and an initial weight may play a role in the final success
of the vertex. Naturally, there are several ways to model how the
weight influences the success of a vertex. In the model considered
in this paper, individuals arrive in the network with a different
initial number of contacts (given to them at birth) and these
initial numbers form the basis for their future success. Later on,
we shall also discuss other ways of how this effect can be
modeled.

The aim of the current paper is to formulate and analyze a model
that combines the rich-get-richer and the rich-by-birth effect.
The model is a preferential attachment model where the number of
edges added upon the addition of a new vertex is a
random variable associated to the vertex. This indeed gives the
desired combination of preferential attachment and
vertex-dependency of degrees upon vertex-birth. For bounded
initial degrees, the model is
included in the very general class of preferential attachment
models treated in \cite{cf}, but the novelty of the model lies in
that the initial degrees can have an arbitrary distribution. In
particular, we can take the weight distribution to be a power law,
which gives a model with two ``competing" power laws: the power
law caused by the preferential attachment mechanism and the power
law of the initial degrees. In such a situation it is indeed not clear
which of the power laws will dominate in the resulting
degrees of the graph. Our main result implies that the most
heavy-tailed
power law wins, that is, the degrees in the resulting graph will
follow a power law with the same exponent as the initial degrees in case
this is smaller than the exponent induced by the preferential
attachment, and with an exponent determined by the preferential
attachment in case this is smaller.

The proof of our main result requires finite moment of order
$1+\varepsilon$ for the initial degrees. However, we believe that the
conclusion is true also in the infinite mean case. More
specifically, we conjecture that, when the distribution of the
initial degrees is a power law with infinite mean, the degree sequence in
the graph will obey a power law with the same exponent as the the
one of the initial degrees. Indeed, the power law of the
initial degrees will
always be the ``strongest" in this case, since preferential
attachment mechanisms only seem to be able to produce power laws
with finite mean. In reality, power laws with infinite mean are
not uncommon, see e.g.\ Table 3.1 in \cite{NewMan} for some
examples, and hence it is desirable to find a model that can
capture this. Unfortunately, we have not been able to give
a full proof for the infinite mean case, but we present partial
results in Section \ref{sec-12}.

We now proceed with a formal definition of the model and the
formulation of the main results.

\subsection{Definition of the model}
\label{sec-11} The model that we consider is described by a graph
process $\{G(t)\}_{t\geq 1}$. To define it, let $\{W_i\}_{i\geq
1}$ be an i.i.d.\ sequence of positive integer-valued random
variables and let $G(1)$ be a graph consisting of two vertices
$v_0$ and $v_1$ with $W_1$ edges joining them. For $t\geq 2$, the
graph $G(t)$ is constructed from $G(t-1)$ in such a way that a
vertex $v_t$, with associated weight $W_t$, is added to the graph
$G(t-1)$, and the edge set is updated by adding $W_t$ edges
between the vertex $v_t$ and the vertices $\{v_0,v_1,
\ldots,v_{t-1}\}$. Thus, $W_t$ is the {\it initial degree} of
vertex $v_t$. Write $d_0(s),\ldots,d_{t-1}(s)$ for the degrees of
the vertices $\{v_0,v_1, \ldots,v_{t-1}\}$ at time $s\geq t-1$.
The endpoints of the $W_t$ edges emanating from vertex $v_t$ are
chosen independently (with replacement) from
$\{v_0,\ldots,v_{t-1}\}$, and the probability that $v_i$ is chosen
as the endpoint of a fixed edge is equal to
    \eq \label{eq.proportion} \frac{d_i(t-1)+\delta}{\sum_{j=0}^{t-1}
    (d_j(t-1)+\delta)}=\frac{d_i(t-1)+\delta}{2L_{t-1}+t\delta}, \qquad
    0 \leq i \leq t-1, \en
where $L_t=\sum_{i=1}^tW_i$, and $\delta$ is a fixed
parameter of the model. Write $S_{\sss {\rm W}}$ for the support of the
distribution of the initial degrees. To ensure that the above expression
defines a probability, we require that
    \eqn{
    \lbeq{delta-restr}
    \delta+\min\{x: x\in S_{\sss {\rm W}}
    \}>0.
    }
This model will be referred to as the PARID-model (Preferential
Attachment with Random Initial Degrees). Note that, when $W_i\equiv 1$ and
$\delta=0$, we retrieve the original preferential attachment model from
Barab\'asi-Albert \cite{barabasi-1999-286}.

\begin{rem}
{\rm In the PARID-model, we assume that the different edges
of a vertex are attached
in an {\it independent} way, and we also take a simple choice
for the initial graph $G(1)$. However, the proofs given below
are rather insensitive to the precise model definitions, and
can be applied to slightly different settings as well. For example,
the proof can also be applied to the model used in \cite{BRST01}.}
\end{rem}

\begin{rem}
{\rm We shall give special attention to the case where
$\prob(W_i=m)=1$ for some integer $m\geq 1$. This model is closest in
spirit to the Barab\'asi-Albert model, and it turns out that
sharper bounds are possible for the error terms in this case.
These results will be used in \cite{EskHofHoo07}, where we study
the diameter in preferential attachment models.}
\end{rem}

\subsection{Heuristics and main result}
\label{sec-12}
Our main result concerns the degree sequence in the resulting
graph $G(t)$. To formulate it, let $N_k(t)$ be the number of
vertices with degree $k$ in $G(t)$ and define
$p_k(t)=N_k(t)/(t+1)$ as the fraction of vertices with degree $k$.
We are interested in the limiting distribution of $p_k(t)$ as
$t\to\infty$. This distribution arises as the solution of a
certain recurrence relation, of which we will now give a short
heuristic derivation. First note that, obviously,
    \begin{equation}\label{vv_rek}
    \expec [N_k(t)|G(t-1)]=N_k(t-1)+
    \expec[N_k(t)-N_k(t-1)|G(t-1)].
    \end{equation}
Asymptotically, for $t$ large,  it is very  unlikely that a vertex will
be hit by more than one of the $W_t$ edges added upon the addition
of vertex $v_t$. Let us hence ignore this possibility for the
moment. The difference $N_k(t)-N_k(t-1)$ between the number of
vertices with degree $k$ at time $t$ and time $t-1$ respectively,
is then obtained as follows:

\begin{itemize}
\item[(a)] Vertices with degree $k$ in $G(t-1)$ that are hit by
one of the $W_t$ edges emanating from $v_t$ are subtracted
from $N_k(t-1)$. The
conditional probability that a fixed edge is attached to a vertex with degree
$k$ is $(k+\delta)N_k(t-1)/(2L_{t-1}+t\delta)$,
so that the expected number of edges
connected to vertices with degree $k$ is
$W_t(k+\delta)N_k(t-1)/(2L_{t-1}+t\delta)$. This coincides
with the mean number of vertices with degree $k$ in $G(t-1)$
hit by edges from $v_t$, apart from the case when two edges
are attached to the same vertex.

\item[(b)] Vertices with degree $k-1$ in $G(t-1)$ that are hit by
one of the $W_t$ edges emanating from vertex $v_t$ are added to
$N_k(t-1)$. By
reasoning as above, it follows that the mean number of such
vertices is $W_t(k-1+\delta)N_{k-1}(t-1)/(2L_{t-1}+t\delta)$.

\item[(c)] The new vertex $v_t$ should be added if it has degree
$k$, that is, if $W_t=k$.
\end{itemize}

\noindent Combining this gives
    \begin{eqnarray}\label{change}
    &&\expec\left[N_k(t)-N_k(t-1)|G(t-1)\right]\nn\\
    &&\qquad\approx
    \frac{(k-1+\delta)W_t}{2L_{t-1}+t\delta}N_{k-1}(t-1)-\frac{(k+\delta)W_t}
    {2L_{t-1}+t\delta}N_k(t-1)+\mathbf{1}_{\{W_t=k\}},
    \end{eqnarray}
where the approximation sign refers to the fact that we have
ignored the possibility of two or more edges of an arriving vertex
being connected to the same end vertex. Now assume that $p_k(t)$
converges to some limit $p_k$ as $t\to\infty$, so that hence
$N_k(t)\sim (t+1)p_k$. Also assume that the distribution of the
initial degrees has finite mean $\mu$, so that $L_{t-1}/t\to\mu$. Finally,
let $\{r_k\}_{k\ge 1}$ be the probabilities associated with the weight
distribution, that is,
    \eqn{
    \lbeq{rk-def}
    r_k=\prob(W_1=k), \qquad k\geq 1.
    }
Substituting (\ref{change}) into (\ref{vv_rek}) and replacing
$L_{t-1}$ by $\mu t$, we arrive, after taking double expectations,
at
    \eqn{
    \lbeq{approxrec}
    \expec\left[N_k(t)\right]\approx \expec[N_k(t-1)]+
    \frac{(k-1+\delta)}{t \theta}\expec[N_{k-1}(t-1)]
    -\frac{(k+\delta)} {t\theta}\expec[N_k(t-1)]+r_k,
    }
where $\theta=2+\delta/\mu$. Sending $t\to\infty$, and observing that $\expec[N_k(t)]-\expec[N_k(t-1)]\rightarrow p_k$
implies that $\frac 1t\expec[N_k(t)]\rightarrow p_k$,
for all $k$, then yields the recursion
    \begin{equation}\label{rek}
    p_k=\frac{k-1+\delta}{\theta}p_{k-1}-\frac{k+\delta}{\theta}p_k+r_k.
    \end{equation}
By iteration, it can be
seen that this recursion is solved by
    \begin{equation}\label{pk}
    p_k=\frac{\theta}{k+\delta+\theta}\sum_{i=0}^{k-1}r_{k-i}
    \prod_{j=1}^i \frac{(k-j+\delta)}{(k-j+\delta+\theta)},\quad\qquad k\geq 1,
    \end{equation}
where the empty product, arising when $i=0$, is defined to
be equal to one.
Furthermore, since $\{p_k\}_{k\ge 1}$ satisfies (\ref{rek}), we have that
$\sum_{k=1}^\infty p_k=\sum_{k=1}^\infty r_k=1$. Hence, $\{p_k\}_{k\ge 1}$
defines a probability distribution, and the above reasoning
indicates that the limiting degree distribution in the PARID-model
should be given by $\{p_k\}_{k\geq 1}$. Our main result confirms this
heuristics:

\begin{theorem}\label{PARWds}
If the initial degrees $\{W_i\}_{i \geq 1}$ have finite moment of order
$1+\varepsilon$ for some $\varepsilon>0$, then there exists a
constant $\gamma\in(0,\frac 12)$ such that
$$
\lim_{t\toinf}\prob\left(\max_{k\geq 1}|p_k(t)-p_k|\geq t^{-\gamma}\right)=0,
$$
where $\{p_k\}_{k\ge 1}$ is defined in (\ref{pk}). When $r_m=1$ for some
integer $m\geq 1$, then $t^{-\gamma}$ can be replaced by $C\sqrt{\frac{\log{t}}{t}}$ for some sufficiently large constant
$C$.
\end{theorem}

\noindent To analyze the distribution $\{p_k\}_{k\geq 1}$, first
consider the case when the initial degrees are almost surely constant,
that is, when $r_m=1$ for some positive integer $m$. Then $r_j=0$
for all $j\neq m$, and (\ref{pk}) reduces to
$$
p_k=\left\{
\begin{array}{ll} \frac{\theta\Gamma(k+\delta)\Gamma(m+\delta+\theta)}{\Gamma(m+\delta)
\Gamma(k+1+\delta+\theta)} & \mbox{for $k\geq m$};\\
                      0 & \mbox{for $k<m$},
                    \end{array}
            \right.
$$
where $\Gamma(\cdot)$ denotes the gamma-function. By
Stirling's formula, we have that $\Gamma(s+a)/\Gamma(s)\sim s^{a}$
as $s\to\infty$, and from this it follows that $p_k\sim
c k^{-(1+\theta)}$ for some constant $c>0$.
Hence, the degree sequence obeys a power law
with exponent $1+\theta=3+\delta/m$. Note that, by choosing
$\delta>-m$ appropriately, any value of the exponent larger than 2
can be obtained. For other choices of $\{r_k\}_{k\geq 1}$, the behavior of
$\{p_k\}_{k\geq 1}$ is less transparent. The following proposition asserts
that, if $\{r_k\}_{k\geq 1}$ is a power law, then $\{p_k\}_{k\geq 1}$ is a power law
as well. It also gives the aforementioned characterization of the
exponent as the minimum of the exponent of the $r_k$'s and an
exponent induced by the preferential attachment mechanism.

\begin{prop}\label{prop:pk}
Assume that $r_k=\prob(W_1=k)=k^{-\tau_{\sss{\rm W}}}L(k)$ for some
$\tau_{\sss{\rm W}}>2$ and some function $k\mapsto L(k)$ which is slowly
varying. Then $p_k=k^{-\tau}\hat L(k)$ for some slowly varying
function $k\mapsto \hat L(k)$ and with power-law exponent $\tau$ given by
    \begin{equation}\label{tau}
    \tau=\min\{\tau_{\sss{\rm W}},\tau_{\sss {\rm P}}\},
    \end{equation}
where $\tau_{\sss {\rm P}}$ is the power-law exponent of
the pure preferential attachment model given by
$\tau_{\sss {\rm P}}=3+\delta/\mu$.\\
When $r_k$ decays faster than a power law, then (\ref{tau})
remains true with the convention that $\tau_{\sss{\rm W}}=\infty$.
\end{prop}

\noindent In deriving the recursion (\ref{rek}) we assumed that
the initial degrees $\{W_i\}_{i\geq 1}$ have finite mean $\mu$. Assume now that the
mean of the initial degrees is infinite. More specifically, suppose
that $\{r_k\}_{k\geq 1}$
is a power law with exponent $\tau_{\sss{\rm W}}\in[1,2]$. Then,
we conjecture that the main result above remains true.

\begin{conj}
\label{conj-infmean} When $\{r_k\}_{k\geq 1}$ is a power law distribution
with exponent $\tau_{\sss{\rm W}}\in(1,2)$, then the degree
sequence in PARID-model obeys a power law with the same exponent
$\tau_{\sss{\rm W}}$.
\end{conj}

Unfortunately, we cannot quite prove Conjecture
\ref{conj-infmean}. However, we shall prove a slightly weaker
version of it. To this end, write $N_{\geq k}(t)$ for the number of
vertices with degree larger than or equal to $k$ at time $t$, that
is, $N_{\geq k}(t)=\sum_{i=0}^t {\bf 1}_{\{d_i(t)\ge k\}}$, and let $p_{\geq
k}(t)=N_{\geq k}(t)/(t+1)$. Since $d_i(t)\geq W_i$, obviously
    \begin{equation}\label{p_>k_lower}
    \expec[p_{\geq k}(t)]=\frac{\expec[N_{\ge
    k}(t)]}{t+1}\geq\frac{\expec[\sum_{i=1}^t {\bf 1}_{\{W_i\ge k\}}]}{t+1}
    =\prob(W_1\geq k)\frac{t}{t+1}=\prob(W_1\geq k)(1+o(1)),
    \end{equation}
that is, the expected degree sequence in the PARID-model
is always bounded from below by the weight distribution.
In order to prove a related upper bound, we start by
investigating the expectation of the degrees.

\begin{theorem}\label{thm-infmean}
Suppose that $\sum_{k>x} r_k=\prob(W_1>x) = x^{1-\tau_{\sss{\rm
W}}}L(x)$, where $\tau_{\sss{\rm W}}\in(1,2)$ and $x\mapsto L(x)$
is a slowly varying function at infinity. Then, for every
$s<\tau_{\sss{\rm W}}-1$, there exists a constant $C>0$ and a
slowly varying function $x\mapsto l(x)$ such that, for $i\in \{0,
\ldots, t\}$, we have that
    \[
    \expec[d_i(t)^s]\leq C\Big(\frac{t}{i\vee 1}\Big)^{s/(\tau_{\sss \rm W}-1)} \Big(\frac{l(t)}{l(i)}\Big)^s,
    \]
where $x\vee y=\max\{x,y\}$.
\end{theorem}

\noindent Theorem \ref{thm-infmean} gives an upper
bound for the expected degree sequence:

\begin{cor}
\label{cor-infmean}
If $\sum_{k>x} r_k=\prob(W_1>x) = x^{1-\tau_{\sss{\rm W}}}L(x)$, where
$\tau_{\sss{\rm W}}\in(1,2)$ and $x\mapsto L(x)$ is
a slowly varying function at infinity, then, for every $s<\tau_{\sss{\rm W}}-1$, there exists an
$M$ (independent of $t$) such that
    $$
    \expec[p_{\geq k}(t)]\leq M k^{-s}.
    $$
\end{cor}

\noindent \textbf{Proof.} For $\vep>0$ and $s<\tau_{\sss {\rm W}}-1$, it
follows from Theorem \ref{thm-infmean} and Markov's inequality that
    \eqan{
    \expec[p_{\geq k}(t)]&=\frac1{t+1}\sum_{i=0}^t \prob(d_i(t)\ge k)
    = \frac1{t+1}\sum_{i=0}^t \prob(d_i(t)^s \ge k^s)\nn\\
    &\le \frac1{t+1}\sum_{i=0}^t k^{-s}\expec[d_i(t)^s]
    \le k^{-s} \frac{C}{t+1}\sum_{i=0}^t
    \Big(\frac{t}{i\vee 1}\Big)^{s/(\tau_{\sss \rm W}-1)}
    \Big(\frac{l(t)}{l(i)}\Big)^s\leq Mk^{-s},
    }
since, for $s<\tau-1$ and using \cite[Theorem 2, p. 283]{Fell71},
there exists a constant $c>0$ such that
    \[
    \sum_{i=0}^t (i\vee 1)^{-s/(\tau_{\sss \rm W}-1)}l(i)^{-s}
    =c t^{1-\frac{s}{\tau_{\sss \rm W}-1}}l(t)^{-s}(1+o(1)).
    \]
\qed

Combining Corollary \ref{cor-infmean} with (\ref{p_>k_lower})
yields that, when the weight distribution is a power law with
exponent $\tau_{\sss {\rm W}}\in(1,2)$, the only possible power law for
the degrees has exponent equal to $\tau_{\sss{\rm W}}$. This
statement is obviously not as strong as Theorem \ref{PARWds}, but
it does offer convincing evidence for Conjecture
\ref{conj-infmean}. Theorem \ref{thm-infmean} is proved in Section
\ref{sec-3}.

\subsection{Related work}

Before proceeding with the proofs, we describe some related work.
In Section \ref{sec-rel}, the proof of Theorem \ref{PARWds} is
compared to related proofs that have appeared in the literature,
and we refer there for the extensive literature on power laws
for preferential attachment models.
In this section, we describe work on related models.

As mentioned in the introduction, the paper by Cooper and Frieze
\cite{cf} deals with a very general class of preferential
attachment models, including the PARID-model with bounded initial degrees.
Another way of introducing the rich-by-birth effect in a
preferential attachment model, is the \emph{fitness model},
formulated by Barab\'asi and Bianconi
\cite{bianconi-2000-Bose-Einstein, barabasi-2001-Competition}, and
later generalized by Erg\"un and Rodgers \cite{ergun-2002-303}. We
will shortly describe the model of Erg\"un and Rodgers and some
(non-rigorous) results for the degree sequence.

The idea with the model is that vertices have different ability --
referred to as \emph{fitness} -- to compete for edges. More
precisely, each vertex has two types of fitness, a multiplicative
and an additive fitness associated to it. These are given by
independent copies of random variables $\eta$ and $\zeta$,
respectively. The dynamics of the model is then very similar to
the dynamics of the PARID-model. However, instead of adding a
random number of edges together with each new vertex, new vertices
come with a {\it fixed} number $m$ of edges. Also, instead of connecting
an edge to a given vertex with a probability proportional to the
degree plus $\delta$, the probability of connecting to a given
vertex is proportional to the degree times the multiplicative
fitness plus the additive fitness. Thus, the expression
\eqref{eq.proportion} for the probability of choosing $v_i$ ($0
\leq i \leq t-1$) as the endpoint of a fixed edge is replaced by
        \eq \label{eq.proportion'}
        \frac{\eta_i d_i(t-1)+\zeta_i}
        {\sum_{j=0}^{t-1} \eta_jd_j(t-1)+\sum_{j=0}^{t-1}\zeta_j},\,\, 0 \leq i \leq t-1.
        \en
The original fitness model of Barab\'asi and Bianconi is
obtained when $\zeta\equiv0$. If, in addition, the multiplicative
fitness is the same for all vertices, the model reduces further to the
Albert-Barab\'asi model.

The rich-by-birth effect is present since relatively young
vertices, with a small degree, can acquire edges at a high rate if
the multiplicative fitness or the additive fitness is large.
Therefore, this model is sometimes called the {\it fit-get-rich}
model. Excluding trivial choices for the distribution of $\eta$
and $\zeta$, it is not clear before hand if the graph process is
driven by the rich-get-richer effect, by the rich-by-birth effect
or by a combination of them. When the additive fitness is zero,
Barab\'asi and Bianconi \cite{barabasi-2001-Competition} show that
the distribution of the average degree sequence of $G(t)$ depends
on the distribution of $\eta$. For $\eta$ uniformly distributed on
$[0,1]$, they show (non-rigorously) that the degree sequence
$\{p_k\}_{k\geq 1}$ is given by
    $$
    p_k \sim c\frac{k^{-(1+C^*)}}{\log k},
    $$
where $C^*$ is the solution of the equation $\exp\br{-2/C}=1-1/C$
and $c>0$ a constant, that is, the average degree sequence follows a
generalized power law. When $\eta$ is exponentially distributed, numerical
simulations indicate that the degree sequence also behaves like an
exponential distribution. For the general model with non-zero
additive fitness there are no explicit expressions for the
$p_k$'s. See however \cite{ergun-2002-303} for some special cases.
We mention also that \cite{bianconi-2000-Bose-Einstein} provides a
coupling of the fitness model to a so-called Bose gas. This
coupling gives a way of predicting (non-rigorously) whether
the rich-get-richer or the rich-by-birth effect will be
dominant.

\section{Proof of Theorem \ref{PARWds} and Proposition \ref{prop:pk}}
\label{sec.proof.PARWds}

In this section, we prove Theorem \ref{PARWds} and Proposition
\ref{prop:pk}. We start by proving Proposition
\ref{prop:pk}, since the proof of Theorem \ref{PARWds} makes
use of it.

\subsection{Proof of Proposition \ref{prop:pk}}
\label{sec-23}
Recall the definition (\ref{pk}) of $p_k$. Assume that $\{r_k\}_{k\geq 1}$
is a power law distribution with exponent $\tau_{\sss{\rm W}}>2$, that is,
assume that $r_k= L(k) k^{-\tau_{\sss{\rm W}}},$ for some
slowly varying function $k\mapsto L(k)$. We want to show that then $p_k$
is a power law distribution as well, more precisely, we want to
show that $p_k=\hat L(k) k^{-\tau}$, where
$\tau=\min\{\tau_{\sss{\rm W}},1+\theta\}$ and $k\mapsto \hat{L}(k)$
is again a slowly varying function. To this end, first note that the
expression for $p_k$ can be rewritten in terms of
gamma-functions as
\begin{equation}
\label{pkasgamma}
    p_k=\frac{\theta\cdot\Gamma(k+\delta)}{\Gamma(k+\delta+1+\theta)}\sum_{m=1}^{k}
    \frac{\Gamma(m+\delta+\theta)}{\Gamma(m+\delta)}r_{m}.
\end{equation}
By Stirling's formula, we have that
\begin{equation}
\label{abram1}
    \frac{\Gamma(k+\delta)}{\Gamma(k+\delta+1+\theta)}=
    k^{-(1+\theta)}
\left(1+O\big(k^{-1}\big)\right),
\quad k\to \infty,
\end{equation}
and
\begin{equation}
\label{abram2}
    \frac{\Gamma(m+\delta+\theta)}{\Gamma(m+\delta)}=m^\theta
    \left(1+O\big(m^{-1}\big)\right),
\quad m\to \infty.
\end{equation}
Furthermore, by the assumption, $r_{m}=L(m) m^{-\tau_{\sss{\rm
W}}}$. It follows that
\begin{equation}
\label{eq:aspk}
\sum_{m=1}^{k}\frac{\Gamma(m+\delta+\theta)}{\Gamma(m+\delta)}r_{m}
\end{equation}
is convergent as $k\to\infty$ if $\theta-\tau_{\sss{\rm W}}<-1$,
that is, if $\tau_{\sss{\rm W}}>1+\theta$. For such values of
$\tau_{\sss{\rm W}}$, the distribution $p_k$ is hence a power law
with exponent $\tau_{\sss {\rm P}}=1+\theta$. When
$\theta-\tau_{\sss{\rm W}}\geq -1$, that is, when $\tau_{\sss{\rm
W}}\leq \tau_{\sss {\rm P}} $, the series in \eqref{eq:aspk}
diverges and, by \cite[Lemma, p. 280]{Fell71}, it can be seen that
    $$
    k\mapsto \sum_{m=1}^{k}\frac{\Gamma(m+\delta+\theta)}{\Gamma(m+\delta)}r_{m}
    $$
varies regularly with exponent $\theta-\tau_{\sss{\rm W}}+1$.
Combining this with \eqref{abram1} yields that $p_k$ (compare
\eqref{pkasgamma}) varies regularly with exponent $\tau_{\sss{\rm
W}}$, as desired. \hfill$\Box$

\subsection{Proof of Theorem \ref{PARWds}}
The proof of Theorem \ref{PARWds}
consists of two parts: in the first part, we prove that the
degree sequence is concentrated around its mean, and in the second part,
the mean degree sequence is identified. We formulate these results in
two separate propositions -- Proposition \ref{PARWc} and
\ref{PARWe} -- which are proved in Section \ref{sec-a} and \ref{sec-b}
respectively.

The result on the concentration of the degree sequence is as follows:

\begin{prop}\label{PARWc}
If the initial degrees $\{W_i\}_{i\geq 1}$ in the PARID-model have finite moment
of order $1+\varepsilon$, for some $\varepsilon>0$, then there exists a constant $\alpha\in(\frac 12,1)$ such that
    $$
    \lim_{t\toinf} \prob\left(\max_{k\geq 1}\Big|N_k(t)-\mexpec{N_k(t)}\Big|
    \geq t^{\alpha}\right)=0.
    $$
When $r_m=1$ for
some $m\geq 1$, then $t^{\alpha}$ can be replaced by
$C\sqrt{t\log{t}}$ for some sufficiently large $C$.
Identical concentration estimates hold for $N_{\geq k}(t)$.
\end{prop}

\noindent As for the identification of the mean degree sequence,
the following proposition says that the expected number of
vertices with degree $k$ is close to $(t+1)p_k$ for large $t$.
More precisely, it asserts that the difference between $\expec[N_k(t)]$
and $(t+1)p_k$ is bounded, uniformly in $k$, by a constant times
$t^\beta$, for some $\beta\in[0,1)$.

\begin{prop}\label{PARWe}
Assume that the initial degrees $\{W_i\}_{i\geq 1}$ in the PARID-model have
finite moment of order $1+\varepsilon$ for some $\varepsilon>0$,
and let $\{p_k\}_{k \geq 1}$ be defined as in (\ref{pk}). Then
there exist constants $c>0$ and $\beta\in[0,1)$ such that
\begin{equation}
\label{mainprop22}
\max_{k\ge 1}
|\mexpec{N_k(t)}-(t+1)p_k|\leq ct^\beta.
\end{equation}
for all $k$. When $r_m=1$ for
some $m\geq 1$, then the above estimate holds with $\beta=0$.
\end{prop}

With Propositions \ref{PARWc} and \ref{PARWe} at hand it is not
hard to prove Theorem \ref{PARWds}:\medskip

\noindent{\bf Proof of Theorem \ref{PARWds}:}
Combining \eqref{mainprop22} with the triangle inequality, it follows that
    $$
    \prob\left(\max_{k\geq 1}\big|N_k(t)-(t+1)p_k\big|\geq ct^\beta+t^{\alpha}\right)\leq
    \prob\left(\max_{k\geq 1}\big|N_k(t)-\mexpec{N_k(t)}\big|\geq t^{\alpha}\right).
    $$
By Proposition \ref{PARWc}, the right side tends to 0 as $t\to
\infty$ and hence, since $p_k(t)=N_k(t)/(t+1)$, we have that
$$
\lim_{t\toinf}\prob\left(\max_{k\geq 1}|p_k(t)-p_k|\geq
\frac{ct^\beta+t^{\alpha}}{t+1}\right)=0.
$$
The theorem follows from this by picking $0 <\gamma<
1-\max\{\alpha,\beta\}$. Note that, since $0\leq\beta<1$ and
$\frac 12 <\alpha<1$, we have $0<\gamma<\frac 12$. The proof for
$r_m=1$ is analogous. $\hfill\Box$

\subsection{Proof of Proposition \ref{PARWc}}
\label{sec-a}
This proof is an adaption of a martingale argument, which first
appeared in \cite{BRST01}, and has been used for all
proofs of power-law degree sequences since. Here, we correct an incompleteness in the proof of the published paper, as we explain in more detail below. The idea is
to express the difference $N_k(t)-\mexpec{N_k(t)}$ in terms of a Doob
martingale. After bounding the martingale differences, the
Azuma-Hoeffding inequality can be applied to conclude that the
probability of observing large deviations is suitably small,
at least when the initial number of edges has bounded support.
When the initial degrees $\{W_i\}_{i\geq 1}$ are unbounded, and extra coupling
step is required. The argument for $N_{\geq k}(t)$
is identical, so we focus on $N_k(t)$.

We start by giving an argument when $W_i\leq t^{a}$ for all $i\leq
t$ and some $a\in (0,\frac12)$. First note that
    \eqn{
    \lbeq{Nkbd}
    N_k(t) \leq \frac{1}{k}\sum_{l=k}^{\infty} lN_l(t)
    \leq \frac{1}{k}\sum_{l=1}^{\infty} lN_l(t)=\frac{L_t}{k}.
    }
Thus, $\expec[N_k(t)]\leq \mu t/k$. For
$\alpha\in(\frac{1}{2},1)$, let $\eta>0$ be such that
$\eta+\alpha>1$ (the choice of $\alpha$ will be specified in more
detail below). Then, for any $k$ such that $k>t^{\eta}$ for some
$\eta>0$, the event $|N_k(t)-\mexpec{N_k(t)}|\geq t^{\alpha}$
implies that $N_k(t)\geq t^{\alpha}$, and hence that $L_t\geq
kN_k(t)> t^{\eta+\alpha}$. It follows from Boole's inequality that
    \eqn{
    \label{max_bd}
    \prob\left(\max_{k\geq 1}|N_k(t)-\mexpec{N_k(t)}|\geq t^{\alpha}\right)\leq
    \sum_{k=1}^{t^{\eta}}\prob\Big(|N_k(t)-\mexpec{N_k(t)}|\geq
    t^{\alpha}\Big)+\prob(L_t>t^{\eta+\alpha}).
    }
Since $\eta+\alpha>1$ and $L_t/t\to \mu$, the event $L_t>
t^{\eta+\alpha}$ has small probability. To estimate the
probability $\prob\left(|N_k(t)-\mexpec{N_k(t)}|\geq
t^{\alpha}\right)$, for $n\in[t]=\{1, \ldots, t\}$, denote $W_{[n]}=(W_i)_{i=1}^{n}$ and define
	$$
	M_n=\mexpec{N_k(t)\mid G_n, W_{[n]}},
	$$
where $G(0)$ is defined as the empty graph. Here, we explicitly write that we are conditioning on $G(n)$ {\em as well as} $W_{[n]}$, even though $W_{[n]}$ can be retrieved from $G_n$. 

In the published version of this paper, we did not properly take into account that the knowledge that the initial degree of vertex $n$ equals $W_n$ changes all connection probabilities. This is resolved here. First note that, since $\mexpec{M_n}<\infty$, the process is a Doob martingale with respect with $M_t=N_k(t)$ and $M_0=\mexpec{N_k(t)}$, so that
    $$
    N_k(t)-\mexpec{N_k(t)}=M_t-M_0.
    $$
To estimate the increments $M_n-M_{n-1}$, write
	\eqan{
	M_n-M_{n-1}&=\mexpec{N_k(t)\mid G(n), W_{[n]}}-\mexpec{N_k(t)\mid G(n-1), W_{[n-1]}}\nn\\
	&=\mexpec{N_k(t)\mid G(n), W_{[n]}}-\mexpec{N_k(t)\mid G(n-1), W_{[n]}}\label{Doob_split1}\\
	&\qquad +\mexpec{N_k(t)\mid G(n-1), W_{[n]}}-\mexpec{N_k(t)\mid G(n-1), W_{[n-1]}}\label{Doob_split2}.
	}
We treat the terms \eqref{Doob_split1} and \eqref{Doob_split2} separately. In \eqref{Doob_split1}, we condition on $W_{[n]}$, so also on $W_n$, in both expectations, so that we only change the information about the graph. The additional information contained in $G(n)$ compared to $G(n-1)$ consists in how the $W_n$ edges emanating from $v_n$ are attached. This can affect the degrees of at most $2W_n$ vertices. Therefore, we obtain
	\eqn{
	\label{term-correctly-handled}
	\Big|\mexpec{N_k(t)\mid G(n), W_{[n]}}-\mexpec{N_k(t)\mid G(n-1), W_{[n]}}\Big|\leq 2W_n.
	}
The term \eqref{Doob_split2} was omitted in the analysis in the published version of the paper, but we will deal with it now. To this end, let $W_n'$ be a copy of $W_n$, independent of $W_{[t]}$, and write $W_{[n]}'=(W_1, \ldots, W_{n-1}, W_n')$. Then
	$$
	\mexpec{N_k(t)\mid G(n-1), W_{[n-1]}}
	=\mexpec{\mexpec{N_k(t)\mid G(n-1), W_{[n]}'}\mid G(n-1), W_{[n-1]}}.
	$$
Therefore, with $N_k'(t)$ denoting the number of degree $k$ vertices in the model with initial degrees $(W_1, \ldots, W_{n-1}, W_n', W_{n+1}, \ldots, W_t)$ and $d'_{v}(t)$ the degrees in this model, we arrive at	
	\eqan{
	&\mexpec{N_k(t)\mid G(n-1), W_{[n]}}-\mexpec{N_k(t)\mid G(n-1), W_{[n-1]}}\nn\\
	&=\expec\Big[\expec[N_k(t)\mid G(n-1), W_{[n]}]-\expec[N_k'(t)\mid G(n-1), W_{[n]}']\mid G(n-1), W_{[n-1]}\Big].\nn
	}
We conclude that we need to investigate the effect of the change of the initial degree of vertex $n$ from $W_n$ to an independent copy $W_n'$. Note that
	$$
	N_k(t)=\sum_{v=1}^t\indic{d_v(t)=k},
	$$
so that
	\eqan{
	&\expec[N_k(t)\mid G(n-1), W_{[n]}]-\expec[N_k'(t)\mid G(n-1), W_{[n]}']\nn\\
	&\qquad =\sum_{v=1}^t \Big[\prob(d_v(t)=k \mid G(n-1), W_{[n]})-\prob(d_v'(t)=k \mid G(n-1), W_{[n]}')\Big]\nn\\
	&\qquad =\sum_{v=1}^t \Big[\prob(d_v(t)=k \mid G(n-1), W_{[n]}, W_{[n]}')-\prob(d_v'(t)=k \mid G(n-1), W_{[n]}, W_{[n]}')\Big],\label{main-identity}
	}
since the event $\{d_v(t)=k\}$ is independent of $W_n'$, while $\{d_v'(t)=k\}$ is independent of $W_n$, and conditioning on extra independent information does not affect conditional probabilities. We split the sum in \eqref{main-identity} into three contributions, depending on whether $v<n, v=n$ or $v>n$. When $v=n$, we bound
	\eqn{
	\label{contri-1}
	|\prob(d_n(t)=k \mid G(n-1), W_{[n]}, W_{[n]}')-\prob(d_n'(t)=k \mid G(n-1), W_{[n]}, W_{[n]}')|\leq 1.
	}
When $v<n$, we claim that there exists a $C>0$ such that
	\eqn{
	\label{Dv-vs-Dv'-bd}
	\Big|\prob(d_v(t)=k \mid G(t), W_{[t]}, W_{[t]}')-\prob(d_v'(t)=k \mid G(t), W_{[t]}, W_{[t]}')\Big|
	\leq C\frac{d_v(t)+\delta}{t} |W_n-W_n'|,
	}
so that, by taking the conditional expectation given $G(n-1), W_{[n]}, W_{[n]}'$, also
	\eqan{
	\label{Dv-vs-Dv'-bd-b}
	&\Big|\prob(d_v(t)=k \mid G(n-1), W_{[n]}, W_{[n]}')-\prob(d_v'(t)=k \mid n-1G, W_{[n]}, W_{[n]}')\Big|\\
	&\qquad \leq C|W_n-W_n'|  \expec\Big[\frac{d_v(t)+\delta}{t}\mid G(n-1), W_{[n]}, W_{[n]}'\Big].\nn
	}	
Let us explain \eqref{Dv-vs-Dv'-bd} in more detail. The left hand side of \eqref{Dv-vs-Dv'-bd} can be bounded by the (conditional) probability that $d_v'(t)\neq d_v(t)$. The only difference between the degree of $v<n$ in the two models is the additional $|W_n-W_n'|$ edges that are present in one model, but not the other, and $d_v'(t)\neq d_v(t)$ only when one of these edges connects to $v$. We can think of these edges as being paired {\em last}, and then the (conditional) probability of attaching to vertex $v$ is
	\eqn{
	\label{cond-prob-attach}
	\frac{d_v(t)+\delta}{2\tilde{L}_t+\delta t},
	}
where
	$$
	\tilde{L}_t=(W_n\wedge W_n')+\sum_{i\in [t]\setminus \{n\}} W_i.
	$$
By assumption $\delta +\min\{x\colon x\in S_{\sss W}\}>0$  (recall \refeq{delta-restr}), where $S_{\sss W} $ denotes the support of $W$. Hence, the denominator in \eqref{cond-prob-attach} above is $2\tilde{L}_t+\delta t=\tilde{L}_t+(\tilde{L}_t+\delta t)\geq t$, where we also use that $W\geq 1$. This explains \eqref{Dv-vs-Dv'-bd}. 

Summing \eqref{Dv-vs-Dv'-bd-b} over all $v< n$, we arrive at
	\eqan{
	\label{contri-2}
	&\sum_{v=1}^{n-1} \Big|\prob(d_v(t)=k \mid G(n-1), W_{[n]}, W_{[n]}')-\prob(d_v'(t)=k \mid G(n-1), W_{[n]}, W_{[n]}')\Big|\\
	&\qquad \leq C|W_n-W_n'|\sum_{v=1}^{ n-1} \frac{\expec\big[d_v(t)\mid G(n-1), W_{[n]}, W_{[n]}'\big]+\delta}{t}\leq C|W_n-W_n'|.\nn
	}	
\medskip

We are left with dealing with $v>n$. In this case, the only difference between the degree of $v>n$ in the two models arises from the possibly different attachment probabilities due to the additional $|W_n-W_n'|$ edges that are present in one model, but not in the other. Indeed, the (conditional) attachment probability of vertex $v$ at time $s\in (n,t]$ when vertex $n$ has $W_n$ edges equals
	$$
	\frac{d_v(s-1)+\delta}{2L_{s-1}+\delta (s-1)},
	$$
while in the model where vertex $n$ has $W_n'$ edges, it equals
	$$
	\frac{d_v'(s-1)+ \delta}{2L_{s-1}'+\delta (s-1)},
	$$	
where
	$$
	L_{s-1}=\sum_{i=1}^{s-1} W_{i},
	\qquad
	L_{s-1}'=\sum_{i=1}^{s-1} W_{i}'=L_{s-1}+W_n'-W_n.
	$$
We rewrite
	\eqan{
	&\prob(d_v(t)=k \mid G(n-1), W_{[n]}, W_{[n]}')-\prob(d_v'(t)=k \mid G(n-1), W_{[n]}, W_{[n]}')\nn\\
	&\quad=\prob(d_v(t)=k, d_v'(t)\neq k \mid G(n-1), W_{[n]}, W_{[n]}')-\prob(d_v'(t)=k, d_v(t)\neq k \mid G(n-1), W_{[n]}, W_{[n]}').\nn
	}	
Therefore, at least one of the attachments to $v$ needs to have been made for the $G(n-1), W_{[n]}$ model, and not for the $G(n-1), W_{[n]}'$ model, or the other way around. We consider the {\em first} time $s$ where such an attachment was made differently in the two models, so that at that time  $s-1$, we have $d_v(s-1)=d_v'(s-1)\leq k$. 
Consequently, 	
	\eqan{
	&\prob(d_v(t)=k, d_v'(t)\neq k \mid G(n-1), W_{[n]}, W_{[n]}')\\
	&\qquad\leq
	\sum_{s=n+1}^t \expec\Big[(d_v(s-1)+\delta) \Big|\frac{1}{2L_{s-1}+\delta (s-1)}-\frac{1}{2L'_{s-1}+\delta (s-1)}\Big|\mid G(n-1), W_{[n]}, W_{[n]}'\Big]\nn\\
	&\qquad =\sum_{s=n+1}^t \expec\Big[(d_v(s-1)+\delta) \frac{2|L_{s-1}-L_{s-1}'|}{(2L_{s-1}+\delta (s-1))(2L'_{s-1}+\delta (s-1))}\mid G(n-1), W_{[n]}, W_{[n]}'\Big]\nn\\
	&\qquad =2 |W_n-W_n'|  \sum_{s=n+1}^t \expec\Big[\frac{d_v(s-1)+\delta}{(2L_{s-1}+\delta (s-1))(2L'_{s-1}+\delta (s-1))}\mid G(n-1), W_{[n]}, W_{[n]}'\Big].\nn
	}
Since $d_v(s-1)\ge 1$, we have that $d_v(s-1)+\delta \leq (1+\delta)d_v(s-1)$ for $\delta > 0$, while $d_v(s-1)+\delta \leq d_v(s-1)$ for $\delta\leq 0$. Hence, with $C(\delta)=2\big((1+\delta)\indic{\delta>0}+\indic{\delta\le 0}\big)$, and summing out over $v>n$,
	\eqan{
	&\sum_{v>n} \prob(d_v(t)=k, d_v'(t)\neq k \mid G(n-1), W_{[n]}, W_{[n]}')\nn\\
	&\qquad \leq C(\delta)  |W_n-W_n'|  \sum_{s=n+1}^t \sum_{v>n} \expec\Big[\frac{{d_v(s-1)} }{(2L_{s-1}+\delta (s-1))(2L'_{s-1}+\delta (s-1))}\mid G(n-1), W_{[n]}, W_{[n]}'\Big]\nn\\
	&\qquad  \leq 2C(\delta) |W_n-W_n'|  \sum_{s=n+1}^t \expec\Big[
	\frac{L_{s-1}}
	{(2L_{s-1}+\delta (s-1))(2L'_{s-1}+\delta (s-1))}
	\mid G(n-1), W_{[n]}, W_{[n]}'\Big].\nn
}
where in the second inequality we used that $\sum_{v>n} d_v(s-1)\leq 2L_{s-1}$.

Recall that, by assumption, $\delta +\min\{x\colon x\in S_{\sss W}\}>0$, where $S_{\sss W}$ denotes the support of $W$, so that $L_{s-1}+\delta (s-1)>0$ and $L'_{s-1}+\delta (s-1)>0$. Hence,
	$$	
	\frac{L_{s-1}}{2L_{s-1}+\delta (s-1)}=\frac{L_{s-1}}{L_{s-1}+(L_{s-1}+\delta (s-1))}\leq 1,
	$$
and
	$$
	2L'_{s-1}+\delta (s-1)\geq L'_{s-1}\geq s-1.
	$$
Therefore,
 	\eqan{
	\label{contri-3}
	&\sum_{v>n} \prob(d_v(t)=k, d_v'(t)\neq k \mid G(n-1), W_{[n]}, W_{[n]}')\\
	&\qquad\leq 
	2C(\delta) |W_n-W_n'| \sum_{s=n+1}^t \Big(\frac{1}{s-1}\Big)\nn\\
	&\qquad \leq 2C(\delta) |W_n-W_n'| \log(t/n) .\nn
	}
By summing \eqref{contri-1},  \eqref{contri-2} and \eqref{contri-3}, and replacing $2C(\delta)$ by a general $C$, we conclude that
	\eqan{
	&\Big|\expec[N_k(t)\mid G(n-1), W_{[n]}]-\expec[N_k'(t)\mid G(n-1), W_{[n]}']\Big|\nn\\
	&\qquad \leq 1+C |W_n-W_n'|\log(t/n).\nn
	}
When taking the expectation w.r.t.\ $W_n'$, and also using that $W_n\leq t^a$, we finally obtain that
	\eqan{
	\label{term-missing}
	|M_n-M_{n-1}|&=\Big|\mexpec{N_k(t)\mid G(n-1), W_{[n]}}-\mexpec{N_k(t)\mid G(n-1), W_{[n-1]}}\Big|\\
	&\leq 1+C |W_n-W_n'|  \log(t/n) \nn\\
	&\leq 1+Ct^a\log(t/n)\equiv c_n(t).\nn
	}
The Azuma-Hoeffding inequality -- see e.g.\ \cite[Section 12.2]{GriSti01} -- now gives us that
	\eqn{
	\label{fin_bd}
	\prob\big(|N_k(t)-\expec[N_k(t)]|\geq b\big)=\prob\big(|M_t-M_0|\geq b\big)\leq 2 e^{-b^2/[8\sum_{n=1}^t c_n(t)^2]}.
	}
	 We compute that
	\eqn{
	\label{sum-cnt-squared}
	\sum_{n=1}^t c_n(t)^2 \leq C\Big[t+t^{1+2a}\Big]\leq  C t^{1+2a}.
	}	
Combining \eqref{fin_bd} and \eqref{sum-cnt-squared} with \eqref{max_bd}, we end up with the estimate.
    \eqn{
    \lbeq{AHW'}
    \prob\Big(\max_{k\geq 1}|N_k(t)-\mexpec{N_k(t)}|\geq t^{\alpha}\Big) \leq
    2 t^{\eta}\exp\left\{-t^{2\alpha-1-2a}/8\right\}
    +\prob(L_t>t^{\eta+\alpha}).
    }
This is the same bound as in the published version of the paper, where the contribution from \eqref{Doob_split2} was omitted. Since $a<1/2$, the above exponential tends to 0 for $\alpha<1$ satisfying that $\alpha>a+1/2$. When the initial degrees are {\it bounded}, so that $W_s\leq m$ almost surely, the bound in \eqref{term-correctly-handled} becomes $2m$, while the bound in \eqref{term-missing} can be improved by replacing $t^a$ by $m$. This means that the bound in \eqref{sum-cnt-squared} becomes $Ct$. Using this, we can deduce that the probability that $\max_{k\geq 1}|N_k(t)-\mexpec{N_k(t)}|$
exceeds $C\sqrt{t\log{t}}$ is $o(1)$ for some $C>0$ sufficiently
large. 

We conclude that we have proved the statement for graphs
satisfying that $W_i\leq t^a$ for some $a\in (0,1/2)$ and all
$i=1, \ldots, t$. Naturally, this assumption may not be true. When
the initial degrees are bounded, the assumption is true, even with
$t^a$ replaced by $m$, but we are interested in graphs having
initial degrees with finite $(1+\vep)$-moments. We next extend the
proof to this setting by a coupling argument.

Write
    \eqn{
    W_i'= W_i\wedge t^{a},\qquad L_s'=\sum_{i=1}^s W_i',
    }
where $x\wedge y=\min\{x,y\}$. Then, the above argument shows that
the PARID-model with initial degrees $\{W_i'\}_{i=1}^t$ satisfies
the claim in Proposition \ref{PARWc}. Denote the graph process
with initial degrees $\{W_i'\}_{i=1}^t$ by $\{G'(i)\}_{i=1}^t$
and, for $i\leq s$, the degree of vertex $i$ in $G'(s)$ by
$d_i'(s)$. We now present a coupling between $\{G(i)\}_{i=1}^t$
and $\{G'(i)\}_{i=1}^t$ which is such that, with high probability,
the number of edges that differ is bounded by $t^{b}$ for some
$b\in (0,1)$.

Define the attachment probabilities in $\{G(i)\}_{i=1}^t$ and
$\{G'(i)\}_{i=1}^t$ by
    \eqn{
    \lbeq{pscoupling}
    p_i(s)=\frac{d_i(s-1)+\delta}{2L_{s-1}+\delta s},
    \quad
    p'_i(s)=\frac{d'_i(s-1)+\delta}{2L_{s-1}'+\delta s}
    .}
Now, we couple the edge attachments such that the $l^{\rm th}$
edge of vertex $s$ in {\it both} graphs is attached to $i$ with
probability $p_i(s)\wedge p_i'(s)$. Otherwise, the edge is {\it
miscoupled}. We shall give a bound on the expected number of
miscouplings. The number of miscouplings in $G(s)$ and $G'(s)$ is
denoted by $U_s$, and is defined in more detail as follows. We
define $U_0=0$ and explain recursively how to construct $U_s$ from
$U_{s-1}$.

The number of miscouplings is adjusted after each
edge which is connected. We consider the
edges to be {\it directed}, and call a
directed edge pointing towards $i$ an
{\it in-edge for $i$},
and a directed edge pointing away from $i$ an
{\it out-edge for $i$}. For convenience
later on, we regard an edge from $s$ to $i$
as {\it both} an in-edge for $i$ as well
as an out-edge for $s$.
By the above definitions, the number
of in-edges of $i$ at time $s$ is the in-degree
of $i$ at time $s$, and the number of out-edges of $i$
at time $s$ is the out-degree of $i$
at time $s$. If we denote the in- and out-degrees
of vertex $i$ in $G(s)$ by $d_{i, {\rm in}}(s)$ and
$d_{i, {\rm out}}(s)$, then clearly
$d_{i}(s)=d_{i, {\rm in}}(s)+d_{i, {\rm out}}(s).$
The same holds for the in- and out-degrees
$d_{i, {\rm in}}'(s)$ and
$d_{i, {\rm out}}'(s)$ of vertex $i$ in $G'(s)$.

The edges which are
attached from vertex $s$ are numbered $1, \ldots, W_s$.
When $W_s>t^a$, then an edge with a number between $t^a$ and $W_s$
adds 2 to $U_{s-1}$, and we call both the in-edge of
$i$ and the out-edge of $s$ as belonging to the
miscoupled set. The size of the miscoupled set
at any time $s=0, \ldots, t$ will be equal to
$U_s$.

When the edge number is in between $1$ and $W_s'$,
then we add 1 to $U_{s-1}$ precisely when the edge
is connected {\it differently} for $G(s)$ and $G'(s)$.
In this case, we say that the
in-edge of $i$ belongs to the miscoupled
set, but the out-edge of $s$ does not. The miscoupled
set remains unchanged when an edge is attached in the
same way in $G(s)$ and in $G'(s)$.

We next define the {\it weight} of every in- and out-edge
of vertex $i$ at time $s$ to be equal to 1.
The {\it total weight} of a vertex $i$ in $G(s)$
at time $s$ is the sum of weights
of the in- and out-edges of $i$ in $G(s)$ plus $\delta$.
The total weight of a vertex $i$ in $G'(s)$ is defined
in a similar manner.

The probabilities in \refeq{pscoupling} are precisely proportional
to the {\it total} weight of the vertex $i$ at time $s-1$. As a
result, for an out-edge of vertex $s$ with number in between 1 and
$W_s'$, a miscoupling occurs with probability equal to
$U_{s-1}/\TW_{s-1}$, where $\TW_s$ is the total weight of all
vertices at time $s$ (i.e., all weights in $G(s)$ and $G'(s)$
combined) plus $\delta s$. To see this, we can choose an edge with
probability equal to the total weight of the end vertex of the
edge divided by $\TW_s$. If this edge is not in the miscoupled
set, then we are done, and the two (directed) edges in $G(s)$ and
$G'(s)$ are equal to an in-edge in the vertex which is connected
to the chosen (directed) edge, and an out-edge from the vertex
$s$. If the chosen (directed) edge is in the miscoupled set, then
it is an edge either for $G(s-1)$ or for $G'(s-1)$, but not for
both, and it is chosen with the correct (conditional) probability.
The above rule then constructs the in- and out-edges corresponding
to the edge we wish to attach. Say the chosen edge is an edge for
$G(s-1)$, then we choose the edge for $G'(s-1)$ from the edges of
$G'(s-1)$ with probability equal to $p_i'(s-1)$. As a result, we
do not create any miscoupling when the initial edge drawn was not
in the miscoupled set. The probability of a miscoupling at time
$s$ is therefore equal to $U_{s-1}/\TW_{s-1}$. Note that $\TW_{s}
\geq 2L_{s}+\delta (s+1)$.

The following lemma bounds the expected value of $U_t$:

\begin{lemma}\label{lem-misc}
There exists a constant $K$ and a $b\in (0,1)$
such that
    \eqn{
    \lbeq{MCtbd}
    \expec[U_t]\leq Kt^b.
    }
\end{lemma}

\noindent{\bf Proof of Lemma \ref{lem-misc}:}
We prove Lemma \ref{lem-misc} by induction.
The induction hypothesis is that, for all $0\leq s\leq t$,
    \eqn{
    \lbeq{MCsbd}
    \expec[U_s] \leq K \left(\frac{s}{t}\right)\cdot t^b.
    }
The bound in \refeq{MCtbd} follows from the one in
\refeq{MCsbd} by substituting $s=t$.

We now prove \refeq{MCsbd}. For $s=0$, we have $U_0=0$, which
initializes the induction hypothesis. To advance the induction
hypothesis, we note that $U_s$ is equal to
$U_{s-1}+2(W_s-W_s')+R_s$, where $R_s$ is the number of out-edges
for $s$ with number in between 1 and $W_s'$ that are miscoupled.
As a result, we have
    \eqn{
    \expec[U_s]=\expec[U_{s-1}]+2\expec[W_s-W_s']
    +\expec[R_s].
    }
By the fact that for each out-edge for $s$ with number
in between 1 and $W_s'$, a miscoupling
occurs with probability equal to $U_{s-1}/\TW_{s-1}$,
we have that
    \eqn{
    \expec[R_s]=\expec\Big[\expec[R_s|W_s]\Big]
    =\expec\Big[W_s' \frac{U_{s-1}}{\TW_{s-1}}\Big]
    =\expec[W_s']\expec\Big[\frac{U_{s-1}}{\TW_{s-1}}\Big],
    }
where the last equality follows from the independence
of $W_s$ and $(U_{s-1},\TW_{s-1})$. Now, we use that
$\TW_{s-1}\geq 2L_{s-1}+\delta (s-1),$ together
with the fact that $L_s$ is concentrated around its mean,
to conclude that $L_{s-1}\geq (\mu-\vep) (s-1)$ with high probability.
Thus,
    \eqn{
    \expec[R_s]=\expec\Big[\expec[R_s|W_s]\Big]
    \leq \frac{\mu}{(s-1)(2\mu+\delta-\vep)}
    \expec[U_{s-1}]+
    \mu\prob\big(L_{s-1}\leq (\mu-\vep) (s-1)\big).
    }
Using the induction hypothesis, we arrive at
\begin{eqnarray}
\expec[R_s]&\le& \frac{\mu}{(s-1)(2\mu+\delta-\vep)} K \left(\frac{s-1}{t}\right) t^b
+\mu\prob(L_{s-1}\leq (\mu-\vep)(s-1))\nn\\
&\le& Kt^{b-1} \frac{\mu}{(2\mu+\delta-\vep)}
+\mu\prob(L_{s-1}\leq (\mu-\vep)(s-1)).
\end{eqnarray}
We further bound
    \eqn{
    \expec[W_s-W_s']=\expec[(W_s-t^a){\bf 1}_{\{W_s>t^a\}}]
    \leq t^{-a \vep} \expec[W_s^{1+\vep}]\leq C t^{-a\vep}.
    }
Therefore, by taking $b-1=-a\vep$, we get that
\begin{eqnarray}
\expec[U_s]&\le&\expec[U_{s-1}]+2\expec[W_s-W'_s]+\expec[R_s]\nn\\
&\le & K(s-1)t^{b-1}+2Ct^{b-1}+Kt^{b-1} \frac{\mu}{(2\mu+\delta-\vep)}
+\mu\prob(L_{s-1}\leq (\mu-\vep)(s-1))\nn\\
&= & Kt^{b-1}\left\{(s-1)+2C/K+\frac{\mu}{(2\mu+\delta-\vep)}\right\}
+\mu\prob(L_{s-1}\leq (\mu-\vep)(s-1))\nn\\
&\le& Ks t^{b-1},
\end{eqnarray}
by noting that $\prob(L_{s}\leq (\mu-\vep)s)$ is exponentially
small in $s$, for $s\to \infty$, and using that, since
$2\mu+\delta>\mu$, we can take $\vep>0$ so small that
$\mu/(2\mu+\delta-\vep)<1$, and, after this, we can take $K$ so
large that
$$
\frac{\mu}{(2\mu+\delta-\vep)}+\frac{2C}{K}<1.
$$

With these choices, we have advanced the induction hypothesis.
$\hfill\Box$\medskip

We now complete the proof of Proposition \ref{PARWc}.
The Azuma-Hoefding argument proves that $N_k'(t)$,
the number of vertices with degree $k$ in $G'(t)$,
satisfies the bound in Proposition \ref{PARWc}, i.e.,
that (recall \refeq{AHW'})
    \eqn{
    \lbeq{AHW'rep}
    \prob\left(\max_{k\geq 1}|N_k'(t)-\mexpec{N_k'(t)}|\geq t^{\alpha}\right) \leq
    2 t^{\eta}\exp\left\{-t^{2\alpha-1-2a}/8\right\}
    +\prob(L_t'>t^{\eta+\alpha}),
    }
for $\alpha\in(\frac{1}{2},1)$ and $\eta>0$ such that
$\alpha+\eta>1$ and $a\in(0,\frac{1}{2})$. Moreover, we have for
every $k\geq 1$, that
    \eqn{
    \lbeq{mainconclcoupling}
    |N_k(t)-N'_k(t)|\leq U_t,
    }
since every miscoupling can change the degree of at most one
vertex. By \refeq{mainconclcoupling} and \refeq{MCtbd}, there is a
$b\in(0,1)$ such that
    \eqn{
    \lbeq{AHcouplingbd1}
    \Big|\expec[N_k(t)]-\expec[N'_k(t)]\Big|\leq \expec[U_t]
    \leq Kt^b.
    }
Also, by the Markov inequality, \refeq{mainconclcoupling} and
\refeq{MCtbd}, for every $\alpha\in (b,1)$, we have that
    \eqn{
    \lbeq{AHcouplingbd2}
    \prob\Big(\max_{k\geq 1}|N_k(t)-N'_k(t)|>t^{\alpha}\Big)
    \leq \prob\big(U_t>t^{\alpha}\big)
    \leq t^{-\alpha} \expec[U_t]=o(1).
    }
Now fix $\alpha\in(b\vee (a+\frac{1}{2}),1)$, where $x\vee
y=\max\{x,y\}$, and decompose
    \eqn{
    \lbeq{AHcouplingbd3}
    \max_{k\geq 1}\big|N_k(t)-\expec[N_k(t)]\big|
    \leq \max_{k\geq 1}\big|N_k'(t)-\expec[N_k'(t)]\big|
    +\max_{k\geq 1}\big|\expec[N_k(t)]-\expec[N_k'(t)]\big|
    +\max_{k\geq 1}\big|N_k(t)-N_k'(t)\big|.
    }
The first term on the right hand side
is bounded by $t^{\alpha}$ with high probability
by \refeq{AHW'rep}, the second term is, for $t$ sufficiently large
and with probability one, bounded by $t^{\alpha}$
by \refeq{AHcouplingbd1}
while the third term is bounded by $t^{\alpha}$ with high probability
by \refeq{AHcouplingbd2}. This completes the proof.
$\hfill\Box$\medskip

\subsection{Proof of Proposition \ref{PARWe}}
\label{sec-b}
For $k\geq 1$, let
    \eqn{
    \lbeq{Nbar-def}
    \bar{N}_k(t)=\expec[N_k(t)|\{W_i\}_{i=1}^t]
    }
denote the expected number of
vertices with degree $k$ at time $t$ given the initial degrees
$W_1,\ldots,W_t,$ and define
    \eqn{
    \lbeq{vep-def}
    \varepsilon_k(t)=\bar{N}_k(t)-(t+1)p_k, \qquad k\geq 1.
    }
Also, for $Q=\{Q_k\}_{k\geq1}$ a sequence of real numbers, define
the supremum norm of $Q$ as $\|Q\|=\sup_{k \geq 1}|Q_k|$. Using
this notation, since $\expec[\bar{N}_k(t)]=\expec[N_k(t)]$, we
have to show that there are constants $c>0$ and $\beta\in[0,1)$
such that

\begin{equation}\label{PARWee}
\|\mexpec{\varepsilon(t)}\|=
\sup_{k\ge 1}
|\expec[\bar{N}_k(t)]-(t+1)p_k|
\leq ct^\beta, \qquad\textrm{for}\quad t=0,1,\ldots,
\end{equation}

\noindent where
$\varepsilon(t)=\{\varepsilon_k(t)\}_{k=1}^{\infty}$.
The plan
to do this is to formulate a recursion for $\varepsilon(t)$,
and then use induction in $t$ to establish (\ref{PARWee}). The
recursion for $\varepsilon(t)$ is obtained by combining a
recursion for $\bar{N}(t)=\{\bar{N}_k(t)\}_{k\geq 1}$,
that will be derived below, and the recursion for
$p_k$ in (\ref{rek}). The hard work then is to
bound the error terms in this recursion;
see Lemma \ref{reclem} below.

Let us start by deriving a recursion for $\bar{N}(t)$.
To this end, for a real-valued
sequence $Q=\{Q_k\}_{k\geq 0}$, with $Q_0=0$, introduce
the operator $T_t$, defined as
(compare to \refeq{approxrec})

\begin{equation}\label{T}
(T_tQ)_k=\left(1-\frac{k+\delta}{2L_{t-1}+t\delta}\right)Q_k+
\frac{k-1+\delta}{2L_{t-1}+t\delta}Q_{k-1},\qquad k \geq 1.
\end{equation}

\noindent When applied to $\bar{N}(t-1)$, the operator $T_t$
describes the effect of the addition of a single edge emanating
from the vertex $v_t$, the vertex $v_t$ itself being excluded from
the degree sequence: There are on the average $\bar{N}_{k-1}(t-1)$
vertices with degree $k-1$ at time $t-1$ and a new edge is
connected to such a vertex with probability
$(k-1+\delta)/(2L_{t-1}+t\delta)$. After this connection is made,
the vertex will have degree $k$. Similarly, there are on the
average $\bar{N}_k(t-1)$ vertices with degree $k$ at time $t-1$.
Such a vertex is hit by a new edge with probability
$(k+\delta)/(2L_{t-1}+t\delta)$, and will then have degree $k+1$.
The expected number of vertices with degree $k$ after the addition
of one edge is hence given by the operator in (\ref{T}) applied to
$\bar N(t)$.

Write $T^n_t$ for the $n$-fold application of $T_t$, and
define $T'_t=T_t^{W_t}$. Then $T'_t$ describes the change
in the expected degree sequence $\bar{N}(t)$ when all the $W_t$
edges emanating from vertex $v_t$ have been connected, ignoring
vertex $v_t$ itself. Hence, $\bar{N}(t)$ satisfies
    \begin{equation}\label{N_krec}
    \bar{N}_k(t)=(T'_t\bar{N}(t-1))_k+{\bf 1}_{\{W_t=k\}}, \qquad k \geq 1.
    \end{equation}

Introduce a second operator $S$ on sequences of real numbers
$Q=\{Q_k\}_{k\geq 0}$, with $Q_0=0$, by (compare to (\ref{rek}))
    \begin{equation}\label{S}
    (SQ)_k=\frac{k-1+\delta}{\theta}Q_{k-1}-\frac{k+\delta}{\theta}Q_{k},
    \qquad k\geq 1,
    \end{equation}
where $\theta=2+\delta/\mu$ and $\mu$ is the expectation
of $W_1$.

The recursion (\ref{rek}) is given by
$p_k=(Sp)_k+r_k$, with initial condition
$p_0=0$. It is solved by $p=\{p_k\}_{k\geq 1}$, as defined in
(\ref{pk}). Observe that
    \begin{eqnarray}
    (t+1)p_k & = & tp_k+(Sp)_k+r_k\nonumber\\
    & = & t(T'_tp)_k+r_k-\kappa_k(t)\label{(t+1)pk}, \qquad k \geq 1,
    \end{eqnarray}
where
    \eqn{
    \lbeq{kappa-def}
    \kappa_k(t)=t(T'_tp)_k-(Sp)_k-tp_k.
    }
Combining \refeq{vep-def}, (\ref{N_krec}) and (\ref{(t+1)pk}),
and using the linearity of $T_t'$, it follows that
$\varepsilon(t)=\{\varepsilon_k(t)\}_{k\geq 1}$ satisfies the
recursion
    \eq
    \label{eq.rec}
    \varepsilon_k(t)=(T'_t\varepsilon(t-1))_k+{\bf 1}_{\{W_t=k\}}-r_k+\kappa_k(t),
    \en
indeed,
\begin{eqnarray*}
 \varepsilon_k(t)&=&{\bar N}_k(t)-(t+1)p_k\\
 &=&(T'_t\bar{N}(t-1))_k+{\bf 1}_{\{W_t=k\}}
 -t(T'_tp)_k-r_k+\kappa_k(t)\\
 &=&(T'_t\varepsilon(t-1))_k+{\bf 1}_{\{W_t=k\}}
 -r_k+\kappa_k(t).
\end{eqnarray*}
Now we define $k_t=\eta t,$ where $\eta\in (\mu, 2\mu+\delta)$.
Since, by \refeq{delta-restr}, $\delta > -\min\{x: x\in S_{\sss
{\rm W}}\}\ge -\mu$, the interval $(\mu, 2\mu+\delta)\neq
\varnothing$. Also, by the law of large numbers, $L_t\leq k_t,$ as
$t\to\infty$, with high probability. Further, we define
$\tilde{\varepsilon}_k(t)=\varepsilon_k(t){\bf 1}_{\{k\leq k_t\}}$
and note that, for $k\leq k_t$, the sequence
$\{\tilde{\varepsilon}_k(t)\}_{k\geq 1}$ satisfies
\eq \label{eq.recrep}
    \tilde\varepsilon_k(t)
    ={\bf 1}_{\{k\leq k_t\}} (T'_t\varepsilon(t-1))_k+
    {\bf 1}_{\{W_t=k\}}-r_k+
    \tilde\kappa_k(t),
    \en
where $\tilde\kappa_k(t)=\kappa_k(t){\bf 1}_{\{k\leq k_t\}}$.
It follows from $\mexpec{{\bf 1}_{\{W_t=k\}}}=r_k$ and the triangle inequality that
    \begin{eqnarray}
    \|\mexpec{\varepsilon(t)}\| & \leq &
    \|\mexpec{\varepsilon(t)-\tilde{\varepsilon}(t)}\|
    +\|\mexpec{\tilde{\varepsilon}(t)}\|\nonumber\\
    & \leq & \|\mexpec{\varepsilon(t)-\tilde{\varepsilon}(t)}\|+
    \|\mexpec{{\bf 1}_{(-\infty, k_t]}(\cdot) T'_t\varepsilon(t-1)}\|+
    \|\mexpec{\tilde\kappa(t)}\|,\label{recer}
    \end{eqnarray}
where ${\bf 1}_{(-\infty, k_t]}(k)={\bf 1}_{\{k\leq k_t\}}$.
Inequality (\ref{recer}) is the key ingredient in the proof of
Proposition \ref{PARWe}. We will derive the following bounds for
the terms in (\ref{recer}).
\begin{lemma}\label{reclem}
There are constants $\const{\tilde\varepsilon}$, $\cvepone$, $\cveptwo$ and
$\const{\tilde\kappa}$, independent of $t$, such that for $t$ sufficiently large
and some $\beta \in [0,1)$,
\begin{itemize}
    \item[\rm{(a)}]
    $\|\mexpec{\varepsilon(t)-\tilde{\varepsilon}(t)}\|
    \leq \frac{\const{\tilde\varepsilon}}{t^{1-\beta}}$,

    \item[\rm{(b)}]
    $\|\mexpec{{\bf 1}_{(-\infty, k_t]}(\cdot)T'_t\varepsilon(t-1) }\|
    \leq \left(1-\frac{\cvepone}{t}\right)
    \|\mexpec{\varepsilon(t-1)}\|+\frac{\cveptwo}{t^{1-\beta}}$,

    \item[\rm{(c)}]
    $\|\mexpec{\tilde\kappa(t)}\| \leq
    \frac{\const{\tilde\kappa}}{t^{1-\beta}}$.
\end{itemize}
When $r_m=1$ for some integer $m\geq 1$, then the above
bounds hold with $\beta=0$.
\end{lemma}

\noindent Given these bounds, Proposition \ref{PARWe} is easily
established.

\medskip
\noindent{\bf Proof of Proposition \ref{PARWe}:} Recall
that we want to establish (\ref{PARWee}). We shall prove
this by induction on $t$. Fix $t_0\in {\mathbb N}$.
We start by verifying the
induction hypothesis for $t\leq t_0$, thus initializing
the induction hypothesis. For any $t\leq t_0$, we have

    \eq \label{eq.induc.c}
    \|\mexpec{\varepsilon(t)}\|  \leq
    \sup_{k\geq 1}\mexpec{\bar{N}_k(t)}+(t_0+1)\sup_{k\geq 1}p_k
    \leq  2(t_0+1),
    \en
since there are precisely $t_0+1$ vertices at time $t_0$ and $p_k\leq
1$. This initializes  the induction hypothesis, when
$c$ is so large that $2(t_0+1)\leq ct^\beta_0$.
Next, we advance the induction hypothesis.
Assume that (\ref{PARWee}) holds at time $t-1$ and apply
Lemma \ref{reclem} to (\ref{recer}) to get that
    \begin{eqnarray*}
    \|\mexpec{\varepsilon(t)}\|
    &\leq&
    \|\mexpec{\varepsilon(t)-\tilde{\varepsilon}(t)}\|+
    \|\mexpec{{\bf 1}_{(-\infty, k_t]}(\cdot)T'_t\varepsilon(t-1)}\|+
    \|\mexpec{\tilde{\kappa}(t)}\|
    \\
    &\leq& \frac{\const{\tilde\vep}}{t^{1-\beta}}+
    \left(1-\frac{\cvepone}{t}\right)c(t-1)^\beta+\frac{\cveptwo}{t^{1-\beta}}+
    \frac{\const{\tilde \kappa}}{t^{1-\beta}}\\
    & \leq & ct^\beta-\frac{c\cdot
    \cvepone-(\cveptwo+\const{\tilde\vep}+\const{\tilde \kappa})}{t^{1-\beta}},
    \end{eqnarray*}
as long as $1-\frac{\cvepone}{t}\geq 0$, which is equivalent to
$t\geq \cvepone$. If we then choose $c$ large so that $c\cdot
\cvepone\geq \cveptwo+\const{\vep}+\const{\tilde \kappa}$ and
$c \geq 2(t_0+1)t_0^{-\beta}$ (recall
\eqref{eq.induc.c}) and $t_0\geq \cvepone$, then we have that
$\|\mexpec{\varepsilon(t)}\|\leq ct^\beta$, and (\ref{PARWee})
follows by induction in $t$.\hfill$\Box$\bigskip

It remains to prove Lemma \ref{reclem}. We shall prove
the statements in Lemma \ref{reclem} (a)-(c) one by one,
starting with (a).

\paragraph{Proof of Lemma \ref{reclem}(a):} We have
$\|\mexpec{\varepsilon(t)-\tilde{\varepsilon}(t)}\|
\leq \mexpec{\|\varepsilon(t)-
\tilde{\varepsilon}(t)\|}$, and, using the definition of
$\tilde{\varepsilon}(t)$, we get that
    \begin{eqnarray*}
    \|\varepsilon(t)-\tilde{\varepsilon}(t)\|& = &
    \sup_{k>k_t}\left|\bar{N}_k(t)-(t+1)p_k
    \right|\leq \sup_{k>k_t}\bar{N}_k(t)+(t+1)\sup_{k>k_t}p_k.
    \end{eqnarray*}
The maximal possible degree of a vertex at time $t$ is
$L_t$, implying that $\sup_{k>k_t}\bar{N}_k(t)= 0,$ when $L_t\leq k_t$.
The latter is true almost surely when $r_m=1$ for some integer $m$, when $t$
is sufficiently large, since for $t$ large
$L_t=mt\leq \eta t=k_t$, where $\eta\in (m,2m+\delta)$, by the
fact that $\mu=m$ and $\delta>-m$.
On the other hand, by \refeq{Nkbd}, with $N_k(t)$ replaced by $\bar{N}_k(t)$ we find
$\bar{N}_k(t)\leq \frac{L_t}{k_t}$ for $k\geq k_t$, and
we obtain that
    \eqn{
    \expec[\sup_{k>k_t}\bar{N}_k(t)]\leq (k_t)^{-1}
    \expec[L_t {\bf 1}_{\{L_t>k_t\}}].
    }
With $k_t= \eta t$ for some $\eta\in (\mu, 2\mu+\delta)$, we have that
    \eqn{
    \expec[L_t {\bf 1}_{\{L_t>k_t\}}]\leq
    k_t^{-\vep}\expec[L_t^{1+\vep} {\bf 1}_{\{L_t>k_t\}}]\leq
    k_t^{-\vep}\expec[|L_t-\mu t|^{1+\vep}]
    +(\mu t)^{1+\vep}k_t^{-\vep}\prob(L_t>k_t),
    }
 and, by the Markov inequality
    $$
    \prob(L_t>k_t)\le \prob\Big(|L_t-\mu t|^{1+\vep}>(k_t-\mu t)^{1+\vep}\Big)
    \leq
    (k_t-\mu t)^{-(1+\vep)}\expec[|L_t-\mu t|^{1+\vep}].
    $$
Combining the two latter results, we obtain
    \begin{equation}
    \label{verbMark}
    \expec[L_t {\bf 1}_{\{L_t>k_t\}}]\leq
    k_t^{-\vep}
    \Big(
    1+\Big(\frac{\mu}{\eta-\mu}\Big)^{1+\vep}
    \Big)\expec[|L_t-\mu t|^{1+\vep}].
    \end{equation}

To bound the last expectation, we will use a consequence of the
Marcinkiewicz-Zygmund inequality, see e.g\ \cite[Corollary 8.2 in
\S3]{Gut1}, which runs as follows. Let $q\in[1,2]$, and suppose
that $\{X_i\}_{i\geq 1}$ is an i.i.d.\ sequence with
$\mexpec{X_1}=0$ and $\mexpec{|X_1|^q}<\infty$. Then there exists a
constant $c_q$ depending only on $q$, such that
    \eq \label{Marcinkiewicz-Zygmund}
    \expec\Big[\Big|\Big(\sum_{i=1}^t X_i\Big)^q\Big|\Big] \leq c_q t
    \mexpec{|{X_1}^q|}.
    \en
Applying (\ref{Marcinkiewicz-Zygmund}) with $q=1+\vep$, we
obtain
    \eqn{
    \expec[\sup_{k>k_t}\bar{N}_k(t)]\leq
    k_t^{-(1+\vep)} \Big(
 1+\Big(\frac{\mu}{\eta-\mu}\Big)^{1+\vep}
 \Big)
    \expec[|L_t-\mu t|^{1+\vep}]
    \leq c_{1+\vep} t^{-\vep}.
    }
Furthermore, since by Proposition \ref{prop:pk}, we have $p_k\leq
c k^{-\gamma}$, for some $\gamma>2$ (see also (\ref{tau})), we have that $\sup_{k>k_t}p_k\leq c t^{-\gamma}$
for some constant $c$. It follows that
    \begin{eqnarray*}
    (t+1)\sup_{k>k_t}p_k& \leq &
    \frac{\const{p}}{t^{\gamma-1}},
    \end{eqnarray*}
and, since $\gamma > 2$, part (a) is established with
$\const{\tilde\vep}=c_{1+\vep}+\const{p}$, and $1-\beta=(\vep\wedge\gamma)-1$.
\qed

\paragraph{Proof of Lemma \ref{reclem}(b):}
Moving on to (b), we will start by showing that for $t$ sufficiently large,
    \begin{equation}\label{(b')}
    \|\mexpec{{\bf 1}_{(-\infty, k_t]}(\cdot)
    T_t\varepsilon(t-1) }\|\leq
    \left(1-\frac{\cvepone}{t}\right)
    \|\mexpec{{\bf 1}_{(-\infty, k_t]}(\cdot)\varepsilon(t-1))}\|
    +\frac{\cvepthree}{t^{1-\beta}},
    \end{equation}
which is (b) when we condition on $W_t=1$.  We shall extend
the proof to the case where $W_t\geq 1$ at a later stage.
To prove (\ref{(b')}), we shall
prove a related bound, which also proves useful
in the extension to $W_t\geq 1$. Indeed, we shall prove,
for any real-valued sequence
$Q=\{Q_k\}_{k\geq 0}$, satisfying (i) $Q_0=0$ and (ii)
    \eqn{
    \lbeq{assbd}
    \sup_{k\geq 1} |k+\delta||Q_k|\leq \const{\sss Q}L_{t-1},
    }
that there exists a $\beta\in (0,1)$
(independent of $Q$) and a constant $c>0$ such that for $t$ sufficiently large,
    \eqn{
    \lbeq{bQ}
    \|\mexpec{{\bf 1}_{(-\infty, k_t]}(\cdot)
    T_tQ}\|\leq
    \left(1-\frac{\cvepone}{t}\right)
    \|\mexpec{{\bf 1}_{(-\infty, k_t]}(\cdot)Q}\|
    +\frac{c\const{\sss Q}}{t^{1-\beta}}.
    }
Here we stress that $Q$ can be {\it random},
for example, we shall apply \refeq{bQ} to
$\vep(t-1)$ in order to derive (\ref{(b')}).

In order to prove \refeq{bQ}, we recall
that
    \eqn{
    \expec[
    (T_tQ)_k]
    =\expec\left[
    \left(1-\frac{k+\delta}{2L_{t-1}+t\delta}\right)Q_k+
    \frac{k-1+\delta}{2L_{t-1}+t\delta}Q_{k-1}\right], \qquad k \geq 1.
    }
In bounding this
expectation we will encounter a problem in that
$Q_{k}$, which is allowed to be random,
and $L_{t-1}$ are not independent (for example when
$Q=\vep(t-1)$). To get
around this, we add and subtract the expression on the right hand
side but with the random quantities replaced by their
expectations, that is, for $k\geq 1$, we write
    \begin{eqnarray}
    \mexpec{\br{T_tQ}_k}& = &
    \left(1-\frac{k+\delta}{2\mu(t-1)+t\delta}\right)\mexpec{Q_k}+
    \frac{k-1+\delta}{2\mu(t-1)+t\delta}\mexpec{Q_{k-1}}\label{r0}\\
    & & +\hspace{0.2cm} (k+\delta)\mexpec
    {Q_k
    \frac{2L_{t-1}-2\mu(t-1)}{(2L_{t-1}+t\delta)(2\mu(t-1)+t\delta)}}\label{r1}\\
    & &+ \hspace{0.2cm}(k+\delta-1)\mexpec{
    Q_{k-1}
    \frac{2\mu(t-1)-2L_{t-1}}{(2L_{t-1}+t\delta)(2\mu(t-1)+t\delta)}}\label{r2}.
    \end{eqnarray}

%

\noindent Note that, when $r_m=1$ for some integer $m\ge 1$, then $L_t=\mu t=mt$.
Hence the terms in (\ref{r1}, \ref{r2}) are
both equal to zero, and only (\ref{r0}) contributes. We first deal
with (\ref{r0}).
Observe that $k\le k_t$ implies that  $k\leq (2\mu+\delta)(t-1)$ for $t$ sufficiently large.
Moreover, $k\leq k_t=\eta t$, where $\eta\in (\mu, 2\mu+\delta)$, so that
we have, for $t$ sufficiently large and all $k\le k_t$,
    \eqn{
    \lbeq{posprob}
    1-\frac{k+\delta}{2\mu(t-1)+t\delta}\geq 0.
    }
It follows that, for $t$ sufficiently large,
    \eqan{
    \lbeq{contraction}
    &\sup_{k\leq k_t}
    \left|\left(1-\frac{k+\delta}{2\mu(t-1)+t\delta}\right)\mexpec{Q_k}+
    \frac{k-1+\delta}{2\mu(t-1)+t\delta}\mexpec{Q_{k-1}}\right|\\
    &\qquad \leq
    \left(1-\frac{1}{2\mu(t-1)+t\delta}\right)
    \|\mexpec{{\bf 1}_{(-\infty, k_t]}(\cdot) Q}\|
    \leq \big(1-\frac{\cvepone}{t}\big)\|
    \mexpec{{\bf 1}_{(-\infty, k_t]}(\cdot)Q}\|,\nn
    }
for some constant $\cvepone$. This proves \refeq{bQ} -- with
$\cvepthree=0$ -- when the number of edges is a.s.\ constant
since (\ref{r1}, \ref{r2}) are zero. It remains
to bound the terms (\ref{r1}, \ref{r2}) in
the case where the number of edges is not a.s.\ constant. We will
prove that the supremum over $k$ of the absolute values of both
these terms are bounded by constants divided by $t^{1-\beta}$ for
some $\beta\in[0,1)$. Starting with (\ref{r1}), by using
the assumption (ii) in \refeq{assbd},
as well as $2L_t+\delta t\geq L_t$,
it follows that
    $$
    \sup_{k\geq 1}\left| (k+\delta)\mexpec {Q_k
    \frac{2L_{t-1}-2\mu(t-1)}{(2L_{t-1}+t\delta)(2\mu(t-1)+t\delta)}}\right|\leq
    \frac{c \const{\sss Q}}{t}\mexpec{|L_{t-1}-\mu(t-1)|}.
    $$
To bound the latter expectation, we combine
(\ref{Marcinkiewicz-Zygmund}) for $q=1+\vep$,
with H\"{o}lders inequality, to obtain
    \begin{eqnarray}
    \label{MHbd}
    \mexpec{|L_t-\mu t|} & \leq & \br{\mexpec{\brc||{L_t-\mu
    t}^{1+\eps}}}^{1/(1+\eps)}\leq \br{C_{1+\vep}t \mexpec{\brc||{W_1-\mu }^{1+\eps}}}^{1/(1+\eps)}
    \leq ct ^{1/(1+\eps)},
    \end{eqnarray}
since ${W_i}$ have finite moment of order
$1+\varepsilon$ by assumption, where, without loss of generality,
we can assume that
$\varepsilon\leq 1$. Hence, we have shown that the supremum over
$k$ of the absolute value of (\ref{r1}) is bounded from
above by a constant divided by $t^{1-\beta}$, where
$\beta=1/(1+\varepsilon)$. That the same is true for the term
(\ref{r2}) can be seen analogously. This completes the proof of
\refeq{bQ}.

To prove (\ref{(b')}), we note that, by convention,
$\vep_0(t-1)=0$, so that we only need to prove that
$\sup_{k\geq 1}|k+\delta||\varepsilon_k(t-1)|\leq cL_{t-1}$.
For this, note from \refeq{Nkbd}, the bound $p_k\le ck^{-\gamma},\, \gamma>2,$
and from the lower bound $L_t\ge t$ that
    \begin{eqnarray}
    \sup_{k\geq 1}|k+\delta||\varepsilon_k(t-1)| & \leq & \sum_{k\geq
    1}(k+|\delta|)
    |\varepsilon_k(t-1)|\nn\\
    & \leq & \sum_{k\geq 1}(k+|\delta|) {\bar N}_k(t-1)+t\sum_{k\geq 1}(k+|\delta|) p_k\nn\\
    & \leq & L_{t-1}+|\delta|(t-1)+ct\leq cL_{t-1}\label{kepsbd},
    \end{eqnarray}
for some constant $c$ (in what follows, to avoid tedious
notation, $c$ will denote a positive constant, independent of $k$
and $t$, whose value may vary from line to line). This completes the
proof of (\ref{(b')}).

To complete the proof of Lemma \ref{reclem}(b), we
first show that \refeq{bQ} implies, for every
$1\leq n\leq t$, and all $k \ge 1$,
    \begin{eqnarray}
    \lbeq{extn1}
    \mexpec{{\bf 1}_{\{k\leq k_t\}}
    \big(T_t^n\varepsilon(t-1)\big)_k} &     \leq &
    \Big(1-\frac{\cvepone}{t}\Big)
    \|\mexpec{{\bf 1}_{(-\infty, k_t]}(\cdot)\varepsilon(t-1)}\|
    +\frac{n \cvepthree}{t^{1-\beta}}.
    \end{eqnarray}
To see \refeq{extn1}, we use induction on $n$.
We note that \refeq{extn1} for $n=1$
is precisely equal to (\ref{(b')}), and this initializes
the induction hypothesis. To advance the induction
hypothesis, we note that
    \eqn{
    {\bf 1}_{\{k\leq k_t\}}
    \big(T_t^n\varepsilon(t-1)\big)_k
    ={\bf 1}_{\{k\leq k_t\}}
    T_t \big(Q(n-1)\big)_k,
    }
where $Q_k(n-1)={\bf 1}_{\{k\leq k_t\}} \Big(T_t^{n-1}\varepsilon(t-1)\Big)_k$.
We wish to use \refeq{bQ}, and we first check
the assumptions (i-ii). By definition, $Q_0(n-1)=0$, which establishes
(i). For assumption (ii), we need to do some more work. According to
\eqref{T}, and using that $2L_{t-1}+t\delta> L_{t-1}\ge t-1$,
for $t$ sufficiently large,
    $$
    \sum_{k=1}^\infty (k+|\delta|)(T_tQ)_k\leq
    \left(
    1+\frac1{t}\right)
    \sum_{k=1}^\infty (k+|\delta|)Q_k,
    $$
and hence, by induction,
    $$
    \sum_{k=1}^\infty (k+|\delta|)(T^{n-1}_tQ)_k\leq
    \left(
    1+\frac1{t}\right)^{n-1}
    \sum_{k=1}^\infty (k+|\delta|)Q_k.
    $$
Substituting $Q_k=\vep_k(t-1)$ and using $|\vep_k(t-1)|\leq
N_k(t-1)+tp_k$, yields
    \begin{eqnarray}
    &&\sum_{k\le k_t} (k+|\delta|)(T^{n-1}_tN(t-1))_k
    +
    \sum_{k\le k_t} (k+|\delta|)(T^{n-1}_tp)_k\nn\\
    &&\quad \leq\left(
    1+\frac1{t}\right)^{n-1}\sum_{k=1}^\infty (k+|\delta|)N_k(t-1)
    +
    \left(
    1+\frac1{t}\right)^{n-1}
    \sum_{k=1}^\infty (k+|\delta|)p_k\nn\\
    &&\quad\leq
    \left(
    1+\frac1{t}\right)^{n-1}
    \cdot cL_{t-1},
    \end{eqnarray}
according to (\ref{MHbd}). Using the inequality $1+x\le e^x, \, x\ge 0$, together
with $n\le t$, this in turn yields,
    \begin{eqnarray}
    \sup_{k\geq 1}|k+\delta||Q_k(n-1)| \leq \exp(1)c L_{t-1}
    \lbeq{kQn-1bd},
    \end{eqnarray}
which implies assumption (ii).

By the induction hypothesis, we have that,
for $k\leq k_t$,
    \eqn{
    \mexpec{Q_k(n-1)} \leq
    \left(1-\frac{\cvepone}{t}\right)
    \|\mexpec{{\bf 1}_{(-\infty, k_t]}(\cdot)\varepsilon(t-1)}\|
    +\frac{(n-1)\cvepthree}{t^{1-\beta}},
    }
so that we obtain
    \begin{eqnarray}
    \lbeq{extn}
    \mexpec{{\bf 1}_{\{k\leq k_t\}}
    \big(T_t^n\varepsilon(t-1)\big)_k} \leq
    \left(1-\frac{\cvepone}{t}\right)
    \|\mexpec{{\bf 1}_{(-\infty, k_t]}(\cdot)\varepsilon(t-1)}\|
    +\frac{(n-1) \cvepthree+c\const{\sss 1}}{t^{1-\beta}},
    \end{eqnarray}
which advances the induction hypothesis when $\cvepthree> c \const{\sss 1}$.

By \refeq{extn}, we obtain that, for $W_t\leq t$,
    \begin{eqnarray*}
    \mexpec{{\bf 1}_{\{k\leq k_t\}}
    \big(T'_t\varepsilon(t-1)\big)_k
    \big|W_t} & \leq &
    \left(1-\frac{\cvepone}{t}\right)
    \|\mexpec{\varepsilon(t-1)|W_t}\|+\frac{W_t \cvepthree}{t^{1-\beta}}\\
        &=&\left(1-\frac{\cvepone}{t}\right)\|\mexpec{\varepsilon(t-1)}\|
    +\frac{W_t \cvepthree}{t^{1-\beta}},
    \end{eqnarray*}
where we use that $\varepsilon(t-1)$ is independent of $W_t$.
When $W_t>t$, we bound, using \refeq{kQn-1bd},
    \eqn{
    \lbeq{T'vepbd}
    \sup_{k\leq k_t}|\big(T'_t\varepsilon(t-1)\big)_k|\leq \const{\sss 4}L_t,
    }
so that
    \eqn{
    \mexpec{{\bf 1}_{\{k\leq k_t\}}
    \big(T'_t\varepsilon(t-1)\big)_k\big|W_t}
    \leq \left(1-\frac{\cvepone}{t}\right)\|\mexpec{\varepsilon(t-1)}\|
    +\frac{W_t \cvepthree}{t^{1-\beta}}
    + \const{\sss 4}\expec[L_t {\bf 1}_{\{W_t>t\}}\big|W_t].
    }
The bound in (b) follows from this by taking expectations on both sides,
using
    \eqn{
    \expec[L_t {\bf 1}_{\{W_t>t\}}]=\mu(t-1)\prob(W_t>t)+ \expec[W_t {\bf 1}_{\{W_t>t\}}]
    \leq \Big(\frac{\mu}{t^{\vep}}+\frac1{t^\vep}\Big)\expec[W_t^{1+\vep}],
    }
after which we use that $\beta=1/(1+\vep)\geq 1-\vep$ and
choose the constants appropriately. This completes the proof
of Lemma \ref{reclem}(b).\qed

%
%

\paragraph{Proof of Lemma \ref{reclem}(c):}
For part (c) of the lemma, recall that
    \eqn{
    \tilde\kappa_k(t)=\kappa_k(t){\bf 1}_{\{k\leq k_t\}}
    \qquad \text{with}
    \qquad
    \kappa_k(t)=t((T_t'-I)p)_k-(Sp)_k,
    }
where $T_t$ is defined in (\ref{T}),  $T_t'=T_t^{W_t}$,
$S$ is defined in (\ref{S}), and where $I$ denotes the
identity operator. In what follows, we will assume that $k\leq
k_t$, so that $\tilde\kappa_k(t)=\kappa_k(t)$. We start by proving a
trivial bound on $\kappa_k(t)$. By (\ref{eq.rec}), we have that
\eq 
    \kappa_k(t)=\varepsilon_k(t)-(T'_t\varepsilon(t-1))_k-{\bf 1}_{\{W_t=k\}}+r_k,
    \en
where $\sup_{k\geq 1}|\varepsilon_k(t)|\leq cL_{t}$ by (\ref{kepsbd})
and $\sup_{1\leq k\leq
k_t}|(T'_t\varepsilon(t-1))_k|\leq \const{\sss 4}L_t$ by
\refeq{T'vepbd}, so that hence
    \eqn{
    \lbeq{trivbdkappa}
    \sup_{k\leq k_t} |\kappa_k(t)|\leq \const{\eta} L_t.
    }
(recall that $k_t=\eta t$ where $\eta\in(\mu,2\mu+\delta)$). For
$x\in [0,1]$ and $w\in {\mathbb N}$, we denote
$$
f_k(x;w)=\big((I + x(T_t-I))^w p\big)_k.
$$
Then $\kappa_k(t)=\kappa_k(t;W_t)$, where
    \eqn{
    \lbeq{kappafrewr}
    \kappa_k(t;w)=t[f_k(1;w)-f_k(0;w)]-(Sp)_k,
    }
and $x\mapsto f_k(x;w)$
is a polynomial in $x$ of degree $w$. By a Taylor expansion around $x=0$,
    \eqn{
    \lbeq{fTayl}
    f_k(1;w)=p_k+w\big((T_t-I)p\big)_k+\frac{1}{2} f_k''(x_k;w),
    }
for some $x_k\in (0,1)$, and, since $I + x(T_t-I)$ and $T_t-I$
commute,
    $$
    f_k''(x;w)=w(w-1) \Big((I + x(T_t-I))^{w-2}(T_t-I)^2p\Big)_k.
    $$
We next claim that, on the event $\{k_t \leq 2L_{t-1}+(t-1)\delta\}$,
    $$
    \sup_{k\leq k_t} \Big|\big((I + x(T_t-I))Q\big)_k\Big|\leq \sup_{k\leq k_t}|Q_k|.
    $$
Indeed, $I + x(T_t-I)=(1-x)I+xT_t$, so the claim follows when
$\sup_{k\leq k_t} |(T_tQ)_k|\leq \sup_{k\leq k_t} |Q_k|$.
The latter is the case, since, on the
event that $k+\delta\leq 2L_{t-1}+t\delta$, and arguing
as in \refeq{contraction}, we have
\begin{eqnarray*}
    \sup_{k\leq k_t} |(T_tQ)_k|
    & \leq & \sup_{k\leq k_t}\Big[
    \Big(1-\frac{k+\delta}{2L_{t-1}+t\delta}\Big)|Q_k|+
    \frac{k-1+\delta}{2L_{t-1}+t\delta}|Q_{k-1}|\Big]\\
    & \leq & \Big(1-\frac{1}{2L_{t-1}+t\delta}\Big)
    \sup_{k\leq k_t} |Q_k|.
\end{eqnarray*}

\noindent Since $k\leq k_t$, the inequality $k+\delta\leq
2L_{t-1}+t\delta$ follows when $k_t \leq 2L_{t-1}+(t-1)\delta$.

As a result, on the event $\{k_t \leq 2L_{t-1}+(t-1)\delta\}$, we have that
    \eqa
    \label{f''bd}
    \max_{x\in [0,1]} \sup_{k\leq k_t}
    |f_k''(x;w)|\leq w(w-1) \sup_{k\leq k_t}\big|\big((T_t-I)^2p\big)_k\big|.
    \ena
Now recall the definition (\ref{S}) of the operator $S$, and note
that, for any sequence $Q=\{Q_k\}_{k=1}^\infty$, we can write
    \eqn{
    \lbeq{Tt-Ieq}
    ((T_t-I)Q)_k=\frac{\theta}{(2L_{t-1}+t\delta)}(SQ)_k
    =\frac{1}{t\mu}(SQ)_k+(R_tQ)_k,
    }
where the remainder operator $R_t$ is defined as
    \eqn{
    \lbeq{Rt-def}
    (R_tQ)_k
    =\left(\frac{k+\delta}{2t\mu+t\delta}-\frac{k+\delta}{2L_{t-1}+t\delta}\right)Q_k
    +
    \left(\frac{k-1+\delta}{2L_{t-1}+t\delta}
    -\frac{k-1+\delta}{2t\mu+t\delta}\right)Q_{k-1}.
    }
Combining \refeq{kappafrewr}, \refeq{fTayl},
\eqref{f''bd} and \refeq{Tt-Ieq}, on the
event $\{k_t \leq 2L_{t-1}+(t-1)\delta\}$ and
uniformly for $k\leq k_t$, we obtain that
    \eqn{
    \label{restrkappabd}
    \kappa_k(t;w)\leq \left(\frac{w}{\mu}-1\right)(Sp)_k+
    w t\sup_{k\leq k_t}|(R_tp)_k|+\frac12 w(w-1)t
    \sup_{k\leq k_t}\big|\big((T_t-I)^2p\big)_k\big|,
    }
together with a similar lower bound with minus signs in front
of the last two terms. Indeed,
\begin{eqnarray*}
 \kappa_k(t;w)&=&t[f_k(1;w)-f_k(0;w)]-(Sp)_k\\
 &=&tw\big((T_t-I)p\big)_k+\frac{1}{2} t f_k''(x_k;w)-(Sp)_k\\
 &=&\frac{wt}{\mu t}(Sp)_k+wt(Rp)_k-(Sp)_k+\frac{1}{2} t f_k''(x_k;w),
 \end{eqnarray*}
and \eqref{restrkappabd} follows from this identity and \eqref{f''bd}.

With (\ref{restrkappabd}) at hand, we are now ready to complete the
proof of (c). We start by treating the case where $r_m=1$ for some integer
$m\ge 1$. In this case, with $w=W_t=m=\mu$, we have that
$(\frac{w}{\mu}-1)(Sp)_k\equiv 0$. Furthermore, the inequality
$k_t \leq 2L_{t-1}+(t-1)\delta$ is true almost surely when $t$ is
sufficiently large. Hence, we are done if we can bound the last two
terms in (\ref{restrkappabd}) with $w=W_t$. To do this, note that,
by the definition (\ref{T}) of $T_t$ and the fact that
$2L_{t-1}+t\delta\geq \eta t$ with $\eta>0$,
    \eqn{
    \lbeq{bdT-I}
    \sup_{k\ge 1} \Big|\big((T_t-I)Q\big)_k\Big| \leq \frac{2}{\eta t} \sup_{k\geq 1} (k+|\delta|) |Q_k|.
    }
Applying \refeq{bdT-I} twice yields that
    $$
    \Big|\big((T_t-I)^2p\big)_k\Big|
    \leq  \frac{4}{\eta^2 t^2} \sup_{k\geq 1} (k+|\delta|)^2 p_k,
    $$
and hence, since by Proposition \ref{prop:pk},
$p_k\leq c k^{-\gamma}$ for some $\gamma>2$, there is a
constant $\const{p}$ such that
    \begin{equation}\label{supk2pk}
    \sup_{k\geq 1}(k+|\delta|)^2 p_k\leq \const{p}.
    \end{equation}
Finally, since $L_t=2m t$, we have that
    $$
    \big|(R_tp)_k\big|\leq \frac{2}{m(t-1)t}\sup_{k\geq 1}
    (k+|\delta|)p_k\leq \frac{2\const{p}}{m(t-1)t}.
    $$
Summarizing, we arrive at the statement that there exists
$c_{m,\delta}$ such that
    $$
    \sup_{k\leq k_t} |\kappa_k(t;m)|\leq \frac{c_{m,\delta}}{t},
    $$
which proves the claim in (c) when $r_m=1$, and with $\beta=0$.

We now move to random initial degrees.
For any $a \in (0,1)$, we can split
\begin{equation}
\label{splitkappa}
    \kappa_k(t)=\kappa_k(t){\bf 1}_{\{W_t\leq t^a\}}
    +\kappa_k(t){\bf 1}_{\{W_t>t^a\}}.
\end{equation}
On the event $\{k_t \leq 2L_{t-1}+(t-1)\delta\}$, the first term
of \eqref{splitkappa} can be bounded by the right side of
(\ref{restrkappabd}), i.e.,
    \eqa
    &&\kappa_k(t){\bf 1}_{\{W_t\leq t^a\}}\nn\\
    &&\quad\,\leq \Big(
    (W_t/\mu-1)(Sp)_k+
        tW_t\sup_{k\leq k_t}|(R_tp)_k|+\frac{W_t(W_t-1)}{2}t
        \sup_{k\leq k_t}\big|\big((T_t-I)^2p\big)_k\big| \Big){\bf 1}_{\{W_t\leq t^a\}}\nn,
    \ena
with a similar lower bound where the last two terms have a minus sign.
From \refeq{trivbdkappa}, we obtain the upper bound
    $$
    \kappa_k(t){\bf 1}_{\{W_t>t^a\}}
    \le
    \const{\eta} L_t {\bf 1}_{\{W_t>t^a\}}.
    $$
Combining these two upper bounds with \eqref{splitkappa}, and adding the
term $(W_t/\mu-1)(Sp)_k{\bf 1}_{\{W_t> t^a\}}$ to the right side, yields that
on the event that $\{k_t \leq 2L_{t-1}+(t-1)\delta\}$,
    \eqan{
    \lbeq{kappabdclose1}
    \kappa_k(t) &\leq  \left(\frac{W_t}{\mu}-1\right)(Sp)_k+
    t W_t{\bf 1}_{\{W_t\leq t^a\}}\sup_{k\leq k_t}|(R_tp)_k|\\
     &\qquad+
     t W_t^2{\bf 1}_{\{W_t\leq t^a\}}\sup_{k\leq  k_t}\big|\big((T_t-I)^2p\big)_k\big|
    +{\bf 1}_{\{W_t>t^a\}}C_\eta L_t,\nn
    }
and  similarly we get as a lower bound,
\eqan{
    \lbeq{kappabdclose2}
    \kappa_k(t) &\geq  \left(\frac{W_t}{\mu}-1\right)(Sp)_k-
    t W_t{\bf 1}_{\{W_t\leq t^a\}}\sup_{k\leq k_t}|(R_tp)_k|\\
     &\qquad-
     t W_t^2{\bf 1}_{\{W_t\leq t^a\}}\sup_{k\leq  k_t}\big|\big((T_t-I)^2p\big)_k\big|
    -{\bf 1}_{\{W_t>t^a\}}\left(C_s\Big|\frac{W_t}{\mu}-1\Big|+C_\eta L_t\right),\nn
    }
where we used that $\sup_{k\ge 1} |(Sp)_k|\le C_s$.
We use \refeq{kappabdclose1} and \refeq{kappabdclose2} on  $\{k_t \leq 2L_{t-1}+(t-1)\delta\}$,
and \refeq{trivbdkappa} on the event $\{k_t > 2L_{t-1}+(t-1)\delta\}$
to arrive at
    \eqan{
    \lbeq{kappabdthere}
    \kappa_k(t) &\leq  \left(\frac{W_t}{\mu}-1\right)(Sp)_k+
    t W_t{\bf 1}_{\{W_t\leq t^a\}}\sup_{k\leq k_t}|(R_tp)_k|\\
     &\,+
     t W_t^2{\bf 1}_{\{W_t\leq t^a\}}\sup_{k\leq  k_t}\big|\big((T_t-I)^2p\big)_k\big|
    +\big({\bf 1}_{\{W_t>t^a\}}+{\bf 1}_{\{k_t >2L_{t-1}+(t-1)\delta\}}\big)
    \left(C_s\Big|\frac{W_t}{\mu}-1\Big|+C_\eta L_t\right),\nn
    }
with a similar lower bound where the last three terms have a minus sign.
We now take expectations on both sides of \refeq{kappabdthere} and take advantage of the equality $\mexpec{W_t/\mu}=1$
and the property that $(Sp)_k$ is deterministic, so that the first term
on the right side drops out. Moreover, using that $W_t$ and
$L_{t-1}$ are independent, as well as that $k_t >2L_{t-1}+(t-1)\delta$
implies that $L_{t-1}\leq k_t$, we arrive at
    \begin{eqnarray}
    |\mexpec{\kappa_k(t)}| & \leq &
    \expec\Big[{\bf 1}_{\{W_t>t^a\}}
    \Big(C_s\Big|\frac{W_t}{\mu}-1\Big|+C_\eta t\Big)\Big]\lbeq{kappa1}\\
    &&+\Big(C_\eta
    \big(k_t+\expec[W_t]\big)+C_s\expec\big[\big|\frac{W_t}{\mu}-1\big|\big]\Big)
    \prob(k_t >2L_{t-1}+(t-1)\delta)\lbeq{kappa2}\\
    &&+t\expec\Big[\sup_{k\leq k_t}|(R_tp)_k|\Big]
    \expec\Big[W_t {\bf 1}_{\{W_t>t^a\}}\Big]\lbeq{kappa3}\\
    &&+
    t
    \expec[W_t^2 {\bf 1}_{\{W_t\leq t^a\}}]
    \expec\Big[\sup_{k\leq k_t}\big|\big((T_t-I)^2p\big)_k\big|\Big].
    \lbeq{kappa4}
    \end{eqnarray}
We now bound each of these four
terms one by one. To bound \refeq{kappa1},
we use that $W_t$ has
finite $(1+\varepsilon)$-moment, to obtain that
    $$
    \mexpec{{\bf 1}_{\{W_t>t^a\}}W_t}
    =\mexpec{{\bf 1}_{\{W_t>t^a\}}W_t^{-\vep} W_t^{1+\vep}}
    \leq t^{-a\vep} \mexpec{W_t^{1+\vep}}
    =O(t^{-a\vep}),
    $$
and,
    $$
    t\mexpec{{\bf 1}_{\{W_t>t^a\}}}
    =t\prob\big(W_t^{1+\vep}>t^{a(1+\vep)}\big)
    \leq t^{1-a(1+\vep)}\mexpec{W_t^{1+\vep}}
    =O(t^{1-a(1+\vep)}),
    $$
which bounds \refeq{kappa1} as
    \eqn{
    \label{bounda}
    \expec\Big[{\bf 1}_{\{W_t>t^a\}}\Big(C_s\Big|\frac{W_t}{\mu}-1\Big|+C_\eta t\Big)\Big]
    =O(t^b),
    }
with $b=\max\{-a\vep,1-a(1+\vep)\}$.

To bound \refeq{kappa2}, we use that
when $k_t >2L_{t-1}+(t-1)\delta$, then
$L_{t-1}<\frac12 (\eta t-\delta (t-1))=\frac12 (\eta-\delta)(t-1) +\frac 12 \eta$.
Now, since $\eta\in (\mu, 2\mu+\delta)$, we have that
$\frac12 (\eta-\delta)<\mu$.
Standard Large Deviation theory and the fact that the initial degrees
$W_i$ are non-negative give that the probability
that $L_{t-1}< \alpha (t-1)$, with $\alpha<\mu$, is exponentially small
in $t$. As a result, we obtain that
    \eqn{
    \lbeq{bounda'}
    \Big(C_\eta
    \big(k_t+\expec[W_t]\big)+C_s\expec\big[\big|\frac{W_t}{\mu}-1\big|\big]\Big)
    \prob\Big(k_t >2L_{t-1}+(t-1)\delta\Big)=O(t^{-1}).
    }

To bound \refeq{kappa3}, we use that $2L_{t-1}+t\delta\geq L_{t-1}\geq t-1\geq t/2$, and also use (\ref{supk2pk}), to obtain that
    $$
    \expec\Big[\sup_{k\leq k_t}|(R_tp)_k|\Big]\leq
    \frac{c}{t^2}\expec|L_{t-1}-t\mu|\sup_{k\geq 1} (k+|\delta|)p_k
        \leq \frac{c}{t^2}\expec|L_{t-1}-t\mu|.
    $$
Thus,
    \eqn{
    \label{boundb}
    t\expec\Big[\sup_{k\leq k_t}|(R_tp)_k|\Big]\expec\Big[W_t {\bf 1}_{\{W_t>t^a\}}\Big]\leq
    \frac{c}{t}\expec|L_{t-1}-t\mu|\cdot t^{-a\vep}\leq O\left(t^{-a\vep-\vep/(1+\vep)}\right),
    }
where the final bound follows from (\ref{MHbd}).

Finally, to bound \refeq{kappa4}, note that
    $$
    \expec[W_t^2 {\bf 1}_{\{W_t\leq t^a\}}]
    =\expec[W_t^{1-\vep} W_t^{1+\vep}{\bf 1}_{\{W_t\leq t^a\}}]
    \leq t^{a(1-\vep)} \expec[W_t^{1+\vep}]
    = O\left(t^{a(1-\vep)}\right),
    $$
and, by (\ref{T}) and the fact that $2L_{t-1}+t\delta\geq \eta t$
for some $\eta>0$, we have
    \eqn{
    \mexpec{\sup_{k\leq k_t}\big|\big((T_t-I)^2p\big)_k\big|}
    \leq \frac{c}{t^2}\sup_{k\geq 1} (k+|\delta|)^2 p_k.
    }
This leads to the bound that
    \eqn{
    \label{boundc}
    t \expec[W_t^2 {\bf 1}_{\{W_t\leq t^a\}}]
    \mexpec{\sup_{k\leq k_t}\big|\big((T_t-I)^2p\big)_k\big|}
    \leq O\left(t^{a(1-\vep)-1}\right).
    }
Combining the bounds in (\ref{bounda}), \refeq{bounda'},
(\ref{boundb}) and (\ref{boundc})
completes the proof of part (c) of
Lemma \ref{reclem}, for any $a$ such that $1/(\vep+1)<a<1$.
\hfill$\Box$


\subsection{Discussion and related results}
\label{sec-rel} In this section, we discuss the similarities and
the differences of our proof of the asymptotic degree sequence in
Theorem \ref{PARWds} as compared to other proofs that have
appeared in the literature.

Virtually all proofs of asymptotic power laws in preferential
attachment models use the two steps presented here in Propositions
\ref{PARWc} and \ref{PARWe}. For bounded support of $W_i$,
the concentration result in Proposition \ref{PARWc} and
its proof are identical in all proofs.
For unbounded $W_i$, an additional coupling argument is
required.
The differences arise in the statement and proof of Proposition
\ref{PARWe}. Our Proposition \ref{PARWe} proves a stronger result
than that for $\delta=0$ appearing in \cite{BRST01} for the
fixed number of edges case, and in \cite{cf} in the random number
of edges case, in that the result is valid for a wider range of
$k$ values and the error term is smaller. Indeed, in \cite{cf},
the equivalent of Proposition \ref{PARWe} is proved for $k\leq
t^{1/21}$, and in \cite{BRST01} for $k\leq t^{1/15}$.

In \cite{cf}, also a random number of edges $\{W_i\}_{i\geq 1}$
is allowed. However, it is assumed that the support of
$W_i$ is {\it bounded},
in which case the Azuma-Hoeffding argument presented in
Section \ref{sec-a} simplifies considerably. The nice feature
of allowing an unbounded number of edges is the competition of the
exponents in \eqref{tau}. The model in \cite{cf} is much more
general than the model discussed here, and at every time allows
for the creation of a new vertex with a random number of edges
or the addition of a random number of edges to an old vertex.
Under reasonable assumptions on the parameters of the model,
a power law is proved for the degree sequence of the graph,
indicating that the occurrence of power laws is rather
robust in the model definition (in contrast to the value of
the power-law exponent). Due to the complexity of the model
in \cite{cf},
the asymptotic degree sequence satisfies a more involved
recurrence relation (see \cite[Eq. (2)]{cf})
than the one in \eqref{rek}.

We close this discussion by reviewing some related results.
Similar results for various random graph processes where a {\it
fixed} number of edges is added can be found in \cite{HagWiu06},
where also similar error bounds are proved for models where a
fixed number of edges is added. In \cite{BBCR03}, a {\it directed}
preferential attachment model is investigated, and it is proved
that the degrees obey a power law similar to the one in \cite{BRST01}.
Finally, in \cite{AieChuLu02}, the error bound in
Proposition \ref{PARWc} is proved for $m=1$ for several models.
The result for fixed $\{W_t\}_{t\geq 1}$ and $m>1$ is,
however, not contained there. We intend to make use of this result
in order to study distances in preferential attachment models in
\cite{EskHofHoo07}. For related references, see
\cite{HagWiu06} and \cite{Szym05}. We finally mention the results in
\cite{KatMor05}. There, a scale-free graph process is studied
where, conditionally on $G(t)$, edges are added {\it
independently} with a probability which is proportional to the
degree of the vertex. In this case, as in \cite{BRST01}, the
power-law exponent can only take the value $\tau=3$, but it
can be expected that by incorporating an additive $\delta$-term
as in \eqref{eq.proportion}, the model can be generalized to
$\tau\geq 3$. Since $\delta<0$ is not allowed in this model
(by the independence of the edges, a degree is zero with
positive probability), we expect that only $\tau\geq 3$
is possible.
\section{Proof of Theorem \ref{thm-infmean}}
\label{sec-3}
In this section, we write $F(x)=\prob(W_1\leq x)$,
and assume that $1-F(x)=x^{1-\tau}L(x)$ for some
slowly varying function $x\mapsto L(x)$. Throughout this
section, we write $\tau=\tau_{\sss {\rm W}}$.

From (\ref{eq.proportion}) it is immediate that
    \eqa
    \label{modeldef}
    d_i(t)=d_i(t-1)+X_{i,t},\qquad \textrm{ for } i=0,1,2,\ldots,t-1,
    \ena
where, conditionally on $d_i(t-1)$ and $\{W_j\}_{j=1}^t$, the
distribution of $X_{i,t}$ is binomial with parameters $W_{t}$
and success probability
        \eqn{
        q_i(t)=\frac{d_i(t-1)+\delta}{2L_{t-1}+t\delta}.
        }
Hence, for
$t>i$,
    \begin{align}
    \mexpec{\br{d_i(t)+\delta}^s}
      \nn& =\expec
      [\expec\big[\left(d_i(t-1)+\delta+X_{i,t} \right)^s|d_i(t-1),\{W_j\}_{j=1}^t\big]]
    \\& \leq \mexpec{\br{d_i(t-1)+\delta+\mexpec{X_{i,t}|d_i(t-1),\{W_j\}_{j=1}^t}}^s},
    \end{align}
where we have used the Jensen inequality $
\expec[(a+X)^s]\le (a+\expec[X])^s$, which follows from concavity
of $t\mapsto (a+t)^s$ for $0<s<1$. Next, we substitute
$\mexpec{X_{i,t}|d_i(t-1),\{W_j\}_{j=1}^t}=W_{t}q_i(t)$ and use
the inequality $2L_{t-1}+t\delta \geq L_{t-1}+\delta$, to obtain
that
    \begin{align*}
    \mexpec{\br{d_i(t)+\delta}^s}& \leq
    \mexpec{(d_i(t-1)+\delta)^s\br{1+\frac{W_{t}}{2L_{t-1}+t\delta}}^s}
    \\& \leq     \mexpec{(d_i(t-1)+\delta)^s\br{1+\frac{W_{t}}{L_{t-1}+\delta}}^s}
    =\mexpec{(d_i(t-1)+\delta)^s\br{\frac{L_{t}+\delta}{L_{t-1}+\delta}}^s}.
    \end{align*}
The above two steps can be repeated, for $t>i+1$, to yield
    \eqa
    \mexpec{\br{d_i(t)+\delta}^s}
    &\leq&\nn \mexpec{\mexpec{(d_i(t-1)+\delta)^s\Big|
    d_i(t-2),\{W_j\}_{j=1}^{t}}\br{\frac{L_{t}+\delta}{L_{t-1}+\delta}}^s}
    \\&\nn \leq&
    \mexpec{(d_i(t-2)+\delta)^s\br{\frac{L_{t-1}+\delta}{L_{t-2}+\delta}}^s
    \br{\frac{L_{t}+\delta}{L_{t-1}+\delta}}^s}
    .\nn
    \ena
Thus, by induction, and because $d_i(i)=W_i$, we get that, for all $t>i\ge 1$,
    \eqa
    \label{productbefore}
    \mexpec{\br{d_i(t)+\delta}^s}\le
    \expec\left[ (W_i+\delta)^s \prod_{n=i+1}^t
    \left( \frac{L_n+\delta}{L_{n-1}+\delta} \right)^s
    \right]
    =\expec\left[ (W_i+\delta)^s \br{\frac{L_t+\delta}{L_{i}+\delta}}^s \right].
    \ena
The case $i=0$ can be treated by $\mexpec{\br{d_0(t)+\delta}^s}=\mexpec{\br{d_1(t)+\delta}^s}$,
which is immediate from the definition of $G(1)$.

Define $f(W_i)=(W_i+\delta)^s$  and
    $$
    g(W_i)=\br{\frac{L_t+\delta}{L_{i}+\delta}}^s
    =\br{1+\frac{W_{i+1}+W_{i+2}+\cdots+W_t}{W_1+W_2+\ldots+W_i+\delta}}^s,
    $$
and notice that when we condition on all $W_j,\, 1\le j \le t$,
except $W_i$, then the map $W_i\mapsto f(W_i)$
is increasing in its argument, whereas
$W_i\mapsto g(W_i)$ is decreasing. This implies that,
    \eqa
    \expec[f(W_i)g(W_i)]\le \expec[f(W_i)]\expec[g(W_i)]. \label{negcorprod}
    \ena
Hence,
    \eqa
    \label{productafter}
    \expec[(d_i(t)+\delta)^s]\le
    \expec[ (W_i+\delta)^s]
            \expec\left[ \br{\frac{L_t+\delta}{L_{i}+\delta}}^s  \right]
    \leq \expec[ (W_i+\delta)^s]
            \expec\left[ (L_t+\delta)^s\right] \expec\left[ (L_{i}+\delta)^{-s}  \right],
    \ena
where in the final step we have applied the inequality
\eqref{negcorprod} once more.

For $i\to \infty$,
    \begin{align} \label{misc1}
    \expec\left[ (L_{i}+\delta)^{-s}  \right]&= (1+o(1))\expec\left[L_{i}^{-s}  \right],&
    \expec\left[ (L_{t}+\delta)^{s}  \right]&= (1+o(1))\expec\left[L_{t}^{s}  \right].&
    \end{align}
The moment of order $s$ of $W_i+\delta$ can be bounded by
    \eqa \label{misc2} \mexpec{(W_i+\delta)^s}
    \leq \mexpec{W_i^{s}\br{1+\frac{|\delta|}{W_i}}^s} \leq
    (1+|\delta|)^s\mexpec{W_i^{s}}=(1+|\delta|)^s\mexpec{W_1^{s}},
    \ena
since $W_i\geq 1$.
Combining \eqref{productafter}, \eqref{misc1} and \eqref{misc2} gives
for $i$ sufficiently large and $t>i$,
    \eqn{ \label{bound}
    \expec[(d_i(t)+\delta)^s]\le (1+|\delta|)^s\mexpec{W_1^s}
    \expec\left[L_{i}^{-s}  \right] \expec\left[L_{t}^{s}  \right](1+o(1)).
    }
We will bound each of the terms  $\mexpec{W_1^s}$, $\expec\left[L_{t}^{s}  \right]$
and $\expec\left[L_{i}^{-s}  \right]$ separately.

It is well known that a positive random variable, where the tail
of the survival function is regularly varying with exponent
$1-\tau$ possesses all moments of order $s< \tau-1$,
so that $\mexpec{W_1^s} < \infty$

Take a norming sequence $\{a_n\}_{n\geq 1}$ such that
     \begin{align}
     \label{props}
      a_n=\sup\left\{x: 1-F(x)\geq \frac{1}{n}\right\},
     \end{align}
then it is immediate that $a_n=n^{1/(\tau-1)}l(n)$,
where $n\mapsto l(n)$ is slowly varying. Also, conveniently,
we have that $1-F(a_n)\geq 1/n$.
In the Appendix (cf. Lemma \ref{eltee}, part (a)), we
show that for some constant $C_s$
     \eqn{\label{bound2'}
     \mexpec{L_t^s} = C_s a_t^s(1+o(1)) \leq C_s t^{s/(\tau-1)}l(t)^s.
     }
As a second result, we prove in the Appendix (cf. Lemma \ref{eltee}, part (b)), that,
for $i$ sufficiently large,
    \eqn{
    \label{bound3'}
    \mexpec{L_i^{-s}}\leq C_s a_i^{-s} =C_s i^{-s/(\tau-1)}l(i)^{-s}.
    }
Combining the equations \eqref{bound},
\eqref{bound2'} and \eqref{bound3'}, we obtain
    \eqa
    \label{endresult'}
    \expec[\br{d_i(t)+\delta}^s]\leq
    C\Big(\frac{t}{i\vee 1}\Big)^{s/(\tau-1)}\Big(\frac{l(t)}{l(i)}\Big)^s.
    \ena
Finally, we note that, since $d_i(t)\geq \min\{x: x\in S_{\sss {\rm W}}\}
\equiv \delta+\eta$ where $\eta>0$,
and using \refeq{delta-restr}, we can bound $\expec[d_i(t)^s]\leq (1\vee \eta^{-1})^s\expec[\br{d_i(t)+\delta}^s]$, which
together with \eqref{endresult'} establishes
the proof of Theorem \ref{thm-infmean}
subject to the proof of
Lemma \ref{eltee}, parts (a) and (b).
\qed


\section*{Acknowledgements}
The work of MD, HvdE and RvdH was supported in part by the Netherlands Organisation for Scientific Research (NWO).

\section{Appendix}

\begin{lemma}
\label{eltee}
Let $F(x)=\prob(W_1\leq x)$,
and assume that $1-F(x)=x^{1-\tau}L(x)$, and
let $L_t=\sum_{i=1}^t W_i$ where $\{W_i\}_{i=1}^t$
are i.i.d.\ copies of $W_1$. Then, there exists a
constant $C_s$ such that
$$
({\rm a}) \qquad \mexpec{L_t^s} = C_s a_t ^s(1+o(1)) \leq C_s
t^{s/(\tau_{\sss \rm W}-1)}l(t)^s.
$$
and
$$
({\rm b})\qquad \mexpec{L_t^{-s}}\leq C_s a_t^{-s} =C_s
t^{-s/(\tau_{\sss \rm W}-1)}l(t)^{-s}.
$$
\end{lemma}

\noindent{\bf Proof:} We start with part (a). We bound
$$
L_t^s=(U_t+V_t)^s \leq 2^s (U_t^s+V_t^s),
$$
where \eqa \label{defUV} U_t=\sum_{j=1}^t W_j {\bf
1}_{\{W_j>a_t\}}, \qquad V_t=\sum_{j=1}^t W_j {\bf 1}_{\{W_j\leq
a_t\}} \ena By concavity of $x\mapsto x^s$,
$$
\mexpec{V_t^s}\le(\mexpec{V_t})^s\leq \Big(t
\int_0^{a_t}[1-F(x)]\,dx\Big)^s \sim (2-\tau_{\sss \rm W})^{-s} a_t^s,
$$
according to \cite[Theorem 2(i), p. 448]{Fell71}. For the other
term, we use that
$$
U_t \le X W_{\sss (t)},
$$
where $W_{\sss (t)}=\max_{1\le j\le t} W_j$, the maximum summand
and $X$ the number of $W_j,\, 1\le j\le t$, which are larger than
$a_t$. Then from the H\"older inequality with $p,q>1$ such that
$p^{-1}+q^{-1}=1$,
\begin{equation}
\label{holder} \mexpec{U_t^s}\leq (\mexpec{X^{sp}})^{1/p}
(\mexpec{W_{\sss (t)}^{sq}})^{1/q}.
\end{equation}
From $t[1-F(a_t)]\to 1$, we see that all moments of $X$ are
bounded, so taking $q>1$, but arbitrarily close to $1$, we can assume
that $sq<\tau_{\sss {\rm W}}-1$, and it follows from Theorem 2.1 of
\cite{pickands} that
    $$
    a_t^{-sq}\mexpec{W_{\sss (t)}^{sq}}\to \mexpec{\zeta^{sq}},
    $$
where $\zeta$ is the limit in distribution of $W_{\sss (t)}$,
since $\mexpec{\zeta^{sq}}$ is finite when $sq<\tau_{\sss {\rm W}}-1$. Hence,
\begin{equation}
\label{bovv} \mexpec{U_t^s}\leq (\mexpec{X^{sp}})^{1/p}
(\mexpec{W_{\sss (t)}^{sq}})^{1/q} = O \big((a_t)^{sq}\big)^{1/q}
=O \big((a_t)^{s}\big).
\end{equation}
Combining these bounds yields the bound (a) of the lemma.

We now turn to part (b) of the proof. We use that $L_{t}\geq
W_{\sss (t)}$, so that
    \eqn{
    \expec\left[L_{t}^{-s}\right]\leq \expec\left[W_{\sss (t)}^{-s}\right]
    =-\expec\left[(-Y_{\sss (t)})^{s}\right],
    }
where $Y_i=-W_i^{-1}$ and $Y_{\sss (t)}=\max_{1\leq i\leq t} Y_i$.
Clearly, $Y_i\in [-1,0]$, so that $\expec[(-Y_{1})^{s}]<\infty$.
Also, $a_t Y_{\sss (t)}=-a_t/W_{\sss(t)}$ converges in
distribution to $-E^{-1/(\tau_{\sss {\rm W}}-1)}$, where $E$ is
exponential, so again it follows from Theorem 2.1 of
\cite{pickands} that \eqn{
    \expec\left[(a_t/L_{t})^{s}\right]\leq -\expec\left[(-a_tY_{\sss (t)})^{s}\right]
    \rightarrow \expec[E^{-1/(\tau_{\sss {\rm W}}-1)}]<\infty.
    }

\qed

\end{document}